\documentclass[preprint,3p]{elsarticle}
\journal{J. Comput. Phys.}

\usepackage{stmaryrd}
\usepackage{amsmath,amssymb}
\usepackage{tabularx,epsfig,color,tabularx}
\usepackage{epstopdf}
\usepackage{algorithm,algorithmic,multirow}
\usepackage{graphicx}
\usepackage{caption}
\usepackage{subcaption}
\usepackage{mathtools}
\DeclareMathOperator*{\argmin}{argmin}

\newcommand{\be}{\begin{equation}}
\newcommand{\ee}{\end{equation}}
\newcommand{\ba}{\begin{array}}
\newcommand{\ea}{\end{array}}
\newcommand{\bea}{\begin{eqnarray}}
\newcommand{\eea}{\end{eqnarray}}
\newcommand{\beas}{\begin{eqnarray*}}
\newcommand{\eeas}{\end{eqnarray*}}

\newcommand\bR{\mathbb{R}}
\newcommand{\bx}{{\bf x}}

\newcommand{\bk}{{\bf k}}

\newcommand{\bz}{{\bf 0}}

\newcommand\wrho{\widehat{\rho}}

\newtheorem{remark}{Remark}[section]

\begin{document}

\begin{frontmatter}

\title{Computing the ground state and dynamics of the nonlinear Schr\"{o}dinger equation
with nonlocal interactions via the nonuniform FFT}
\author[nus]{Weizhu Bao}
 \ead{matbaowz@nus.edu.sg}

\address[nus]{Department of Mathematics, National University of
Singapore, Singapore 119076, Singapore}
\ead[url]{http://www.math.nus.edu.sg/\~{}bao/}

\author[njit]{Shidong Jiang}
\ead{shidong.jiang@njit.edu}
\address[njit]{Department of Mathematical Sciences, New
Jersey Institute of Technology, Newark, New Jersey, 07102, USA}

\author[IECL,Inria]{Qinglin Tang}
\ead{qinglin.tang@inria.fr}
%\address[csrc]{Beijing Computational Science Research Center, Beijing 100084, P. R. China}
\address[IECL]{Universit\'e de Lorraine, Institut Elie Cartan de
Lorraine, UMR 7502, Vandoeuvre-l\`es-Nancy, F-54506, France}
\address[Inria]{Inria Nancy Grand-Est/IECL-CORIDA, France}

\author[wpi]{Yong Zhang\corref{5}}
 \ead{yong.zhang@univie.ac.at}
\address[wpi]{Wolfgang Pauli Institute c/o Fak. Mathematik,
University Wien, Oskar-Morgenstern-Platz 1, 1090 Vienna, Austria}

\cortext[5]{Corresponding author.}

%%%%% Begin Abstract %%%%%%%%%%%
\begin{abstract}
We present efficient and accurate numerical methods for computing the ground state and
dynamics of the nonlinear Schr\"{o}dinger equation (NLSE) with nonlocal interactions
based on a fast and accurate evaluation of the long-range
interactions via the nonuniform fast Fourier transform (NUFFT). We begin with a review
of the fast and accurate NUFFT based method in \cite{JGB} for nonlocal interactions
where the singularity of the Fourier symbol of the interaction kernel at the origin can
be canceled by switching to spherical or polar coordinates.
We then extend the method to compute other nonlocal interactions whose Fourier symbols
have stronger singularity at the origin that cannot be canceled by the coordinate transform.
Many of these interactions do not decay at infinity in the physical space, which adds another layer
of complexity since it is more difficult to impose the correct artificial boundary conditions
for the truncated bounded computational domain.
The performance of our method against other existing methods is illustrated numerically,
with particular attention on the effect of the size of the computational domain in the physical space.
Finally, to study the ground state and dynamics of the NLSE, we propose efficient and accurate numerical
methods by combining the NUFFT method for potential evaluation with
 the normalized gradient flow using backward Euler Fourier pseudospectral
discretization and time-splitting Fourier pseudospectral method, respectively.
Extensive numerical comparisons are carried out between these methods and other existing methods
for computing the ground state and dynamics of the NLSE with various nonlocal interactions.
Numerical results show that our scheme
performs much better than those existing methods in terms of both accuracy and efficiency.

\end{abstract}


\end{frontmatter}

\section{Introduction}\setcounter{equation}{0}
In this paper, we present efficient and accurate numerical methods and compare them
with existing numerical methods for computing the ground state and dynamics of
the nonlinear Schr\"{o}dinger equation (NLSE).
In dimensionless form, the NLSE with a nonlocal (long-range) interaction in $d$-dimensions ($d=3,2,1$)
is
\bea \label{sps-eq}
&&i\,\partial_t \psi(\bx,t)= \left[-\frac{1}{2} \Delta  +V(\bx)+ \beta \,\varphi(\bx,t) \right]
\psi(\bx,t),\qquad \bx \in {\mathbb{ R}}^d,\quad t>0,\label{sps}\\
\label{sps-colb}
&&\quad \varphi(\bx,t) = \left(U\ast |\psi|^2\right)(\bx,t), \qquad \qquad \bx \in {\mathbb{ R}}^d, \quad t\ge 0;\label{phi}
\eea
with the initial data
\be \label{sps-ini}
\psi(\bx,t=0) = \psi_0(\bx),\qquad \bx \in {\mathbb{ R}}^d.
\ee
Here, $t$ is time, $\bx$ is the spatial coordinates, $\psi:=\psi(\bx,t)$ is the complex-valued wave-function,
$V(\bx) $ is a given real-valued external
potential, $\beta$ is a dimensionless interaction constant (positive for repulsive interaction and
negative for attractive interaction), and $\varphi:=\varphi(\bx,t)$ is a real-valued nonlocal (long-range)
interaction which is defined as the convolution of an interaction kernel
$U(\bx)$ and the density function $\rho:=\rho(\bx,t)=|\psi(\bx,t)|^2$.
The NLSE  with the nonlocal interaction  (\ref{sps-eq})-(\ref{sps-colb}) has been widely used
in modelling a variety of problems arising from quantum physics and chemistry to
materials science and biology. It is nonlinear, dispersive and time transverse invariant,
i.e., if $V(\bx)\to V(\bx)+\alpha$ and $\varphi(\bx,t)\to \varphi(\bx,t)+\delta$, then
$\psi(\bx,t)\to \psi(\bx,t)e^{-i(\alpha+\delta)t},$ which immediately implies that the physical
observables such as the density $\rho(\bx,t)=|\psi(\bx,t)|^2$ are unchanged. In addition, it conserves
the {\sl mass} and {\sl energy} defined as follows:
\bea\label{norm}
N(\psi(\cdot,t))&:=&\int_{{\mathbb R}^d}|\psi(\bx,t)|^2d\bx\equiv \int_{{\mathbb R}^d}|\psi(\bx,0)|^2d\bx
=\int_{{\mathbb R}^d}|\psi_0(\bx)|^2d\bx=N(\psi_0), \qquad t\ge0,\\
\label{energy}
E(\psi(\cdot,t))&:=&\int_{{\mathbb R}^d} \left[\frac{1}{2}\left|\nabla\psi(\bx,t)
\right|^2 +V(\bx)|\psi(\bx,t) |^2
+\frac{1}{2}\,\beta\,\varphi(\bx,t)|\psi(\bx,t) |^2\right]\,d \bx \equiv E(\psi_0).
\eea

One of the most important nonlocal interactions in applications
is the Coulomb interaction whose interaction kernel in 3D/2D is given as  %\cite{BJNY}
\be\label{Kernal}
 U_{\rm Cou}(\bx)=\left\{\ba{l}
\frac{1}{4\pi\; |\bx|\;}, \\[0.5em]
\frac{1}{2\pi |\bx|}, \\[0.5em]
\ea\right. \qquad \Longleftrightarrow \qquad \widehat{U}_{\rm Cou}(\xi)=
\left\{\ba{ll}
\frac{1}{|\bk|^2}, &\quad d=3,\\[0.5em]
\frac{1}{ |\bk|},  &\quad d=2,\\[0.5em]
\ea\right.\qquad \bx,\bk\in \mathbb{R}^d,
\ee
where $\widehat{f}(\bk) =  \int_{{\mathbb R}^d} f(\bx)\;e^{-i \bk\cdot \bx}\, d\bx$
is the Fourier transform of $f(\bx)$ for $\bx,\bk\in \mathbb{R}^d$.
In 3D,  the Coulomb interaction kernel $U_{\rm Cou}(\bx)$ is
exactly the Green's function of the
Laplace operator and thus the nonlocal Coulomb interaction $\varphi$ in (\ref{sps-colb}) also satisfies the Poisson equation in 3D
\begin{equation}\label{PoiDiff}
-\Delta\, \varphi(\bx,t) = |\psi(\bx,t)|^2,\qquad \bx \in \mathbb R^3, \quad
\qquad \lim_{|\bx|\to\infty}\varphi(\bx,t)=0, \qquad \qquad t\ge0.
\end{equation}
In this case, (\ref{sps-eq})-(\ref{sps-colb}) is also referred as
the 3D Schr\"{o}dinger-Poisson system (SPS)
which was derived from the linear Schr\"{o}dinger equation
for a many-body (e.g., $N$ electrons) quantum system with binary Coulomb interaction between
different electrons via the ``mean field limit'' \cite{BEGMY,BGM,EY}.
It has important applications in modelling semiconductor devices
and calculating electronic structures in materials simulation and design.
On the other hand,   the Coulomb interaction kernel $U(\bx)$ in 2D is the Green's function of the
square-root-Laplace operator instead of the Laplace operator and thus the
nonlocal Coulomb interaction $\varphi$ in (\ref{sps-colb}) also satisfies
the fractional Poisson equation in 2D
\bea
\label{SqrtPoiDiff}
\sqrt{-\Delta} \, \varphi(\bx,t) = |\psi(\bx,t)|^2,\qquad \bx \in \mathbb R^2, \quad
\qquad \lim_{|\bx|\to\infty}\varphi(\bx,t)=0, \qquad \qquad t\ge0.
\eea
In this case, (\ref{sps-eq})-(\ref{sps-colb}) could be obtained from the 3D SPS under an
infinitely strong external confinement in the $z$-direction \cite{BJNY,Rdct3DTo2D}. This model could be
used for modelling 2D materials such as graphene and ``electron sheets'' \cite{HoneycombChenWu}.

Another type of interaction from applications is that the interaction kernel $U(\bx)$ is taken
 as the Green's function of the Laplace operator in 3D/2D/1D \cite{SPMCompare}
\be\label{Kernal21d}
 U_{\rm Lap}(\bx)=\left\{\ba{ll}
\frac{1}{4\pi|\bx|}, &\ d=3,\\[0.5em]
-\frac{1}{2\pi}\ln |\bx|, &\ d=2,\\[0.5em]
-\frac{1}{2}|\bx|, &\ d=1,\\[0.5em]
\ea\right.\quad  \Longleftrightarrow \quad \widehat{U}_{\rm Lap}(\bk)=
\frac{1}{ |\bk|^2},
  \qquad \bx,\bk\in\mathbb{R}^d.
 \ee
When $d=3$, $U_{\rm Lap}(\bx)=U_{\rm Cou}(\bx)$ for $\bx\in{\mathbb R}^3$.
When $d=2$, the nonlocal interaction $\varphi$ in (\ref{sps-colb}) with (\ref{Kernal21d})
satisfies the Poisson equation in 2D with the far-field condition
\begin{equation}\label{PoiDiff3}
-\Delta\, \varphi(\bx,t) =|\psi(\bx,t)|^2,\qquad \bx \in \mathbb R^2,
\qquad \lim_{|\bx|\to\infty}\left[\varphi(\bx,t)+\frac{C_0}{2\pi} \ln|\bx|\right]=0,
\quad t\ge0;
\end{equation}
and when $d=1$ with $\bx=x$, it satisfies the Poisson equation in 1D with the far-field condition
\begin{equation}\label{PoiDiff3ab}
-\partial_{xx}\varphi(x,t) =|\psi(x,t)|^2,\qquad x \in \mathbb R,
\qquad \lim_{x\to\pm\infty}\left[\varphi(x,t)+\frac{1}{2}\left(C_0|x|\mp C_1\right)\right]=0,
\quad t\ge0,
\end{equation}
where  $C_0=\int_{{\mathbb R}^d}|\psi(\bx,t)|^2 d\bx=\widehat{|\psi|^2}({\bf 0},t)\equiv
\int_{{\mathbb R}^d}|\psi_0(\bx)|^2 d\bx=\widehat{|\psi_0|^2}({\bf 0})=N(\psi_0)$
and $C_1=\int_{\mathbb R} x |\psi(x,t)|^2 \,dx=\widehat{(x|\psi|^2)}({\bf 0},t)$,
which indicate  that the nonlocal interaction  $ \varphi(\bx,t)\to -\infty$ as $|\bx| \to\infty$ in 2D/1D.
In fact, when $d=2$ or $d=1,$ (\ref{sps-eq})-(\ref{sps-colb}) with (\ref{Kernal21d}) is also referred as the
2D or 1D SPS.  They could be obtained from the 3D SPS by integrating the 3D Coulomb interaction kernel $U_{\rm Cou}(\bx)$ along  the $z$-line or  $(y,z)$-plane  under the assumption that
 the electrons are uniformly distributed in one or two spatial dimensions, respectively.
The 2D/1D SPS is usually used for modelling 2D ``electron sheets'' and  1D ``quantum wires'',
respectively, as well as lower dimensions semiconductor devices \cite{Mark}.

Recently, the following nonlocal interaction kernels in 2D/1D were obtained from the 3D SPS
under strongly confining external potentials in the $z$-direction and $(y,z)$-plane,
respectively
\be\label{Kernal124d}
 U_{\rm Con}^\varepsilon(\bx)=\left\{\ba{ll}
\frac{2}{(2\pi)^{3/2}} \int_0^\infty
\frac{e^{-\frac{u^2}{2}}}{\sqrt{|\bx|^2+\varepsilon^2u^2}}\,du, & \bx\in \mathbb{R}^2\\[1em]
 \frac{1}{4 } \;\int_0^{\infty}
\frac{e^{-\frac{u}{2}}}{\sqrt{|\bx|^2+\varepsilon^2 u}}\,du, & \bx\in \mathbb{R} \\[1em]
\ea\right. \Longleftrightarrow \widehat{U}_{\rm Con}^\varepsilon(\bk)=
\left\{\ba{ll}
\frac{2}{\pi}\int_{0}^\infty\frac{e^{- \frac{\varepsilon^2s^2}{2}}}
{|\bk|^2+s^2}\,ds, &\bk\in \mathbb{R}^2, \\[1em]
\frac{1}{2}
\int_0^\infty
\frac{e^{-\varepsilon^2s/2}}{|\bk|^2+s}\,ds, &\bk\in\mathbb{R}, \\[1em]
\ea\right.
 \ee
where $0<\varepsilon\ll1$ is a dimensionless constant describing the ratio of the anisotropic
confinement in different directions in the original 3D SPS \cite{BJNY}.
In this case, the convolution (\ref{sps-colb}) for the
nonlocal interaction $\varphi$
can no longer be re-formulated into a partial differential equation. For
other nonlocal interactions considered  in quantum chemistry and dipole Bose-Einstein condensation,
e.g., the dipole-dipole interaction, we refer to \cite{Bao2013,BC7,DipolarBao,JGB} and references therein.

The ground state $\phi_g$ of the NLSE
is  defined as follows:
\be\label{ground}
\phi_g =\argmin_{\phi\in S} E(\phi), \quad \hbox{where}\quad  S:=\{\phi(\bx) \ |\
\|\phi\|^2:=\int_{{\mathbb R}^d} |\phi(\bx)|^2d\bx=1, \ E(\phi)<\infty\}.
\ee
For the existence, uniqueness and exponentially decay properties of
the ground state as well as the
 well-posedness and dynamical properties of the NLSE,
we refer to \cite{Soler,ReviewSPS,BLS,Bao2013,Rdct3DTo2D,Caz,Masaki,Mehats} and references therein.

In order to numerically compute the ground state of (\ref{ground})
and the dynamics of (\ref{sps-eq})-(\ref{sps-colb}), one of the key difficulties
is to efficiently and accurately evaluate the nonlocal interaction (\ref{sps-colb}) with a given density $\rho=|\psi|^2$.
As we know, a natural way to evaluate a convolution is to compute it in the Fourier domain,
i.e., to re-formulate (\ref{sps-colb}) as
\be\label{rform}
\varphi(\bx,t) = \frac{1}{(2\pi)^d}\int_{{\mathbb R}^d} \widehat{U}(\bk)\, \widehat{|\psi|^2}(\bk,t)\, e^{i\bk\cdot\bx}\,
d\bk=\frac{1}{(2\pi)^d}\int_{{\mathbb R}^d} \widehat{U}(\bk)\, \widehat{\rho}(\bk,t)\, e^{i\bk\cdot\bx}\, d\bk,
\qquad \bx\in{\mathbb R}^d, \quad t\ge0.
\ee
And the integral on the right hand side of \eqref{rform} will be truncated on a ractangular box $\Omega$
in ${\mathbb R}^d$, discretized via the trapezoidal rule, and then computed via the fast Fourier transform
(FFT) \cite{BMS}. However,
the accuracy of this approach is hampered by the fact that
the Fourier transform of the interaction kernel
$\widehat{U}(\bk)$ is singular at the origin.
Indeed, for the Coulomb interaction in 3D,
 it is equivalent to solving the Poisson equation (\ref{PoiDiff}) using the Fourier spectral method
on $\Omega$ with periodic boundary conditions.
It is easy to see that this approach introduces an inconsistency
due to the inappropriate  periodic boundary conditions
 as follows:
\be\label{parad}
0<\int_\Omega |\psi(\bx,t)|^2d\bx =-\int_{\Omega}\Delta \varphi(\bx,t)d\bx=-\int_{\partial\Omega}\frac{\partial\varphi}{\partial \rm n}ds=0.
\ee
Thus, this approach suffers from no convergence in terms of the mesh size of partitioning $\Omega$ when $\Omega$ is small and fixed
(a phenomenon known as ``numerical locking" in the literature); and
its convergence is very slow, e.g.,  linearly convergent for the 3D/2D Coulomb interaction,
in terms of the size of $\Omega$ because $\varphi$ decays like $\frac{1}{|\bx|}$.
To overcome this ``numerical locking'',
a numerical method was proposed by imposing
the homogeneous Dirichlet boundary condition on $\partial \Omega$,
and then solving the truncated problem via the
discrete sine transform (DST) \cite{DipJCP,DipolarBao,SPMCompare}. This method avoids numerically the singularity of $\widehat{U}(\bk)$
at the origin $\bk={\bf 0}$ and thus significantly improves the accuracy in the evaluation of the
Coulomb interaction potential. However,
the truncation error of this method still decays only linearly in terms of the size of $\Omega$
due to the slow decaying property of the Coulomb potential. Thus when high accuracy is required,
the bounded computational
domain $\Omega$ must be chosen very large, which increases significantly the computational cost
in both memory and CPU time for evaluating the nonlocal interaction
potential (\ref{sps-colb}) and solving the NLSE (\ref{sps-eq}).
Moreover, for the purpose of solving the NLSE, a much smaller computational domain actually suffices
since the wave-function
$\psi$ decays exponentially fast when $|\bx|\to\infty$ in most applications.
We would also like to point out that this method could not be extended to the cases where
the potential in (\ref{sps-eq}) either does not decay at infinity (for example,
1D/2D cases of \eqref{Kernal21d}) or
cannot be converted to a PDE problem (as in \eqref{Kernal124d}).

Recently, a fast and accurate algorithm was proposed for the evaluation of the Coulomb interaction
(\ref{Kernal}) in 3D/2D via the NUFFT \cite{JGB}. The key observation there is that the singularity
in the Fourier transform of the interaction kernel
$\widehat{U}(\bk)$ at the origin is canceled out with the Jacobian in spherical or polar coordinates,
thus making the integrand in \eqref{rform} smooth. The integral is then approximated
via a high-order quadrature and the resulting discrete summation is evaluated via the NUFFT.
The algorithm has $O(N\log N)$ complexity with $N$ the total number of unknowns in the physical space
and achieves very high accuracy for the evaluation of Coulomb interactions
\cite{JGB}. The main aims of this paper are fourfold:
(i) to extend the algorithm in \cite{JGB} to evaluate the nonlocal interactions
whose Fourier symbols have stronger singularity at the origin which cannot be
canceled by coordinate transform;
(ii) to compare numerically the newly developed NUFFT based method with the
existing numerical methods that are based on either FFT or DST for the evaluation
of these  nonlocal interactions in terms of the size of the computational domain
$\Omega$ and the mesh size of partitioning $\Omega$;
(iii) to propose efficient and accurate numerical methods for computing
the ground state and dynamics of the NLSE with the nonlocal interactions
(\ref{sps-eq})-(\ref{sps-colb}) by
incorporating the algorithm based on the NUFFT for the evaluation of the nonlocal interaction
into the normalized gradient flow method
and the time-splitting Fourier pseudospectral method, respectively,
and (iv) to compare these two new schemes with those existing numerical methods based on
FFT or DST for computing the ground state and dynamics of the NLSE.

The paper is organized as follows.
In Section 2, we briefly review the NUFFT based algorithm in \cite{JGB}
for the evaluation of the Coulomb interaction in 3D/2D, then
extend it to the general nonlocal interaction (\ref{sps-colb}),
including the cases where $U(\bx)$ is taken as either (\ref{Kernal21d})
or (\ref{Kernal124d}).
In Section 3, we present an efficient and accurate numerical method for
computing the ground state of the NLSE (\ref{sps-eq})-(\ref{sps-colb}) by
  coupling the efficient and accurate evaluation of the nonlocal interaction via
the NUFFT and the normalized gradient flow discretized with the backward Euler Fourier pseudospectral
method, and compare the performance of this method and those existing numerical methods.
In Section 4, an efficient and accurate numerical method is proposed for computing
the dynamics of the NLSE by coupling the efficient and accurate evaluation of the nonlocal interaction via
the NUFFT and the time-splitting Fourier pseudospectral method.
Finally, some concluding remarks are drawn in Section 5.

%%%%%%%%%%%%%%%%%%%%%%%%%%%%%%%%%%%%%%%%
%%
%%%%%%%%%%%%%%%%%%%%%%%%%%%%%%%%%%%%%%%%

\section{An algorithm for the evaluation of the nonlocal interaction via the NUFFT}\setcounter{equation}{0}
In this section, we will propose a fast and accurate
evaluation of the nonlocal interaction
\be\label{nonlocal53}
u(\bx)=(U*\rho)(\bx)=\frac{1}{(2\pi)^d}\int_{{\mathbb R}^d} \widehat{U}(\bk)\, \widehat{\rho}(\bk)\,
e^{i\bk\cdot\bx}\,d\bk, \qquad \bx\in {\mathbb R}^d, \qquad d=3,2,1,
\ee
where $\rho:=\rho(\bx)\ge0$ for $\bx\in{\mathbb R}^d$ is a given smooth density function rapidly decaying
at far field and satisfies $C_0:=\widehat{\rho}({\bf 0})=\int_{{\mathbb R}^d}\rho(\bx) d\bx>0$.
We will first briefly review  the algorithm in \cite{JGB} for fast and accurate
evaluation of the Coulomb interactions in 3D and 2D,
and then extend the algorithm to the cases where $U(\bx)$ in (\ref{nonlocal53}) is taken
as either (\ref{Kernal21d}) or (\ref{Kernal124d}).

\subsection{Coulomb interactions in 3D/2D}
When $U(\bx)$ in (\ref{nonlocal53}) is taken as the the Coulomb interaction kernel (\ref{Kernal}), by
truncating the integration domain in  (\ref{nonlocal53}) into a bounded domain and adopting the
spherical/polar coordinates in 3D/2D, respectively, in the Fourier (or phase) space, we have \cite{JGB}
\bea\label{ifft}
u(\bx)&=&\frac{1}{(2\pi)^{d}}\int_{\mathbb{R}^d}  e^{i\,\bk\cdot\bx}\; \widehat{U}_{\rm Cou}(\bk)\;
\widehat{\rho}(\bk)\,d\bk=\frac{1}{(2\pi)^{d}}\int_{\mathbb{R}^d} \frac{1}{|\bk|^{d-1}}\, e^{i\,\bk\cdot\bx}\;
\widehat{\rho}(\bk)\,d\bk \nonumber\\
&\approx&\frac{1}{(2\pi)^d}\int_{|\bk|\leq P}
\frac{1}{|\bk|^{d-1}}\, e^{ i\bk\cdot \bx}\;\wrho(\bk)\,d\bk \nonumber\\
&=&\frac{1}{(2\pi)^d}
\left\{\begin{aligned}
&\int_0^P\int_0^{\pi}\int_0^{2\pi}
e^{ i\bk\cdot \bx} \, \wrho(\bk)\,\sin\theta\; d|\bk| d\theta d\phi, &d=3,\\
&\int_0^P\int_0^{2\pi}
e^{i\bk\cdot \bx} \; \wrho(\bk)\,d|\bk| d\phi, &d=2,\\
\end{aligned}\right.
\qquad \bx\in \Omega\subset{\mathbb R}^d.
\eea
Here, $P=O(1/\varepsilon_0)^{1/n}$, $\varepsilon_0>0$ is the prescribed
precision (e.g., $\varepsilon_0=10^{-10}$), and $n$
is the decaying rate of $\widehat{\rho}(\bk)$ at infinity (i.e., $\widehat{\rho}(\bk)=O(|\bk|^{-n})$ as $|\bk|\to\infty$).
Correspondingly, we choose a bounded domain $\Omega$ large enough such that the truncation error of $\rho(\bx)$
is negligible. It is easy to see that the singularity of the integrand
 at the origin in phase space is removed in spherical or polar coordinates.
Thus, the above integral can be discretized using high order quadratures and the resulting
summation can be evaluated efficiently via the NUFFT.
This leads to an $O(N\log N)+O(M)$ algorithm
where $N$ is the total number of equispaced points in the physical space  and $M$ is the number of nonequispaced points
in the Fourier space. However, although $M$ is roughly the same order as $N$,
the constant in front of $O(M)$ (e.g., $24^d$ for $12$-digit accuracy) is much greater than the constant
in front of $O(N\log N)$. This makes the algorithm considerably slower than the regular FFT,
especially for three dimensional problems.

An improved algorithm is developed to reduce the computational cost in \cite{JGB}.
First, the integral in \eqref{ifft}
is  further  split into two parts via a simple partition of unity:
\bea\label{3.9}
u(\bx)&\approx&\frac{1}{(2\pi)^{d}}\int_{|\bk|\leq P}\frac{1}{|\bk|^{d-1}}\, e^{ i\bk\cdot \bx}\;\wrho(\bk)\,d\bk
\nonumber\\
&=&\frac{1}{(2\pi)^{d}}\int_{|\bk|\leq P}e^{i\bk\cdot \bx}\,\frac{1-p_d(\bk)}{|\bk|^{d-1}}\,\wrho(\bk)\,d\bk
+\frac{1}{(2\pi)^{d}}\int_{|\bk|\leq P}e^{i\bk\cdot \bx}\,\frac{p_d(\bk)}{|\bk|^{d-1}}\,\wrho(\bk)\,d\bk\nonumber\\
&\approx&\frac{1}{(2\pi)^{d}}\int_{\mathcal D}e^{i\bk\cdot \bx}\,w_d(\bk)\,\wrho(\bk)\,d\bk
+\frac{1}{(2\pi)^{d}}\int_{|\bk|\leq P}e^{i\bk\cdot \bx}\,\frac{p_d(\bk)}{|\bk|^{d-1}}\,\wrho(\bk)\,d\bk
:=I_1+I_2,\quad \bx\in \Omega.
\eea
Here, $\mathcal D =\{\bk=(k_1,\ldots, k_d)^T \big| -P\leq k_j \leq P, j = 1,\ldots,d\}$ is a rectangular domain
containing the ball $B$,  the function $p_d(\bk)$ is chosen such that it is a  $C^\infty$ function that decays exponentially fast as $|\bk|\to\infty$  and the function $w_d(\bk):=\frac{1-p_d(\bk)}{|\bk|^{d-1}}$ is smooth for
$\bk\in\bR^d$.

\begin{figure}[!ht]
\centering
\begin{subfigure}[b]{0.35\textwidth}
                \includegraphics[width=\textwidth]{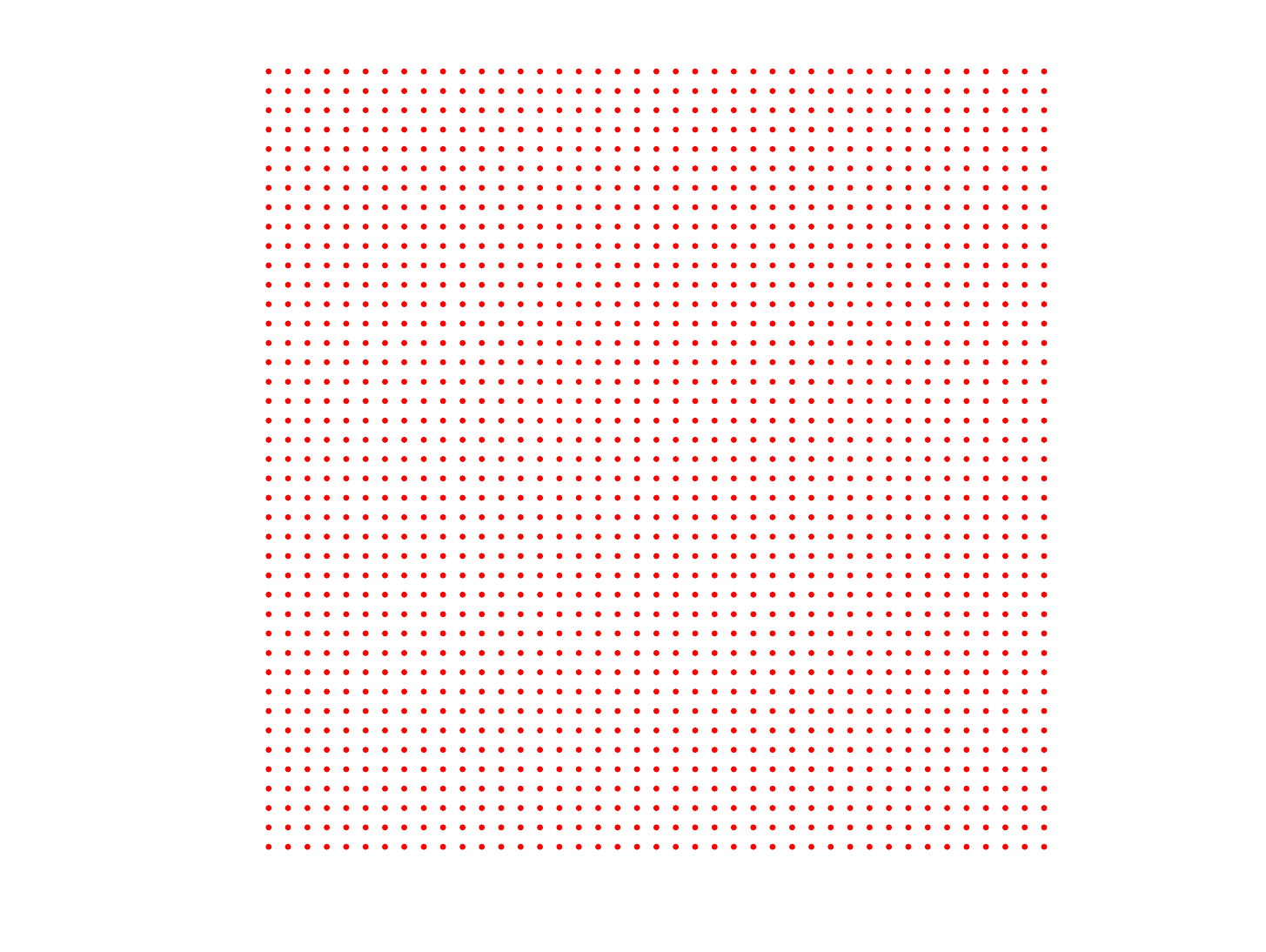}
                \caption{Regular grid}
        \end{subfigure}%
\qquad \qquad
\begin{subfigure}[b]{0.35\textwidth}
                \includegraphics[width=\textwidth]{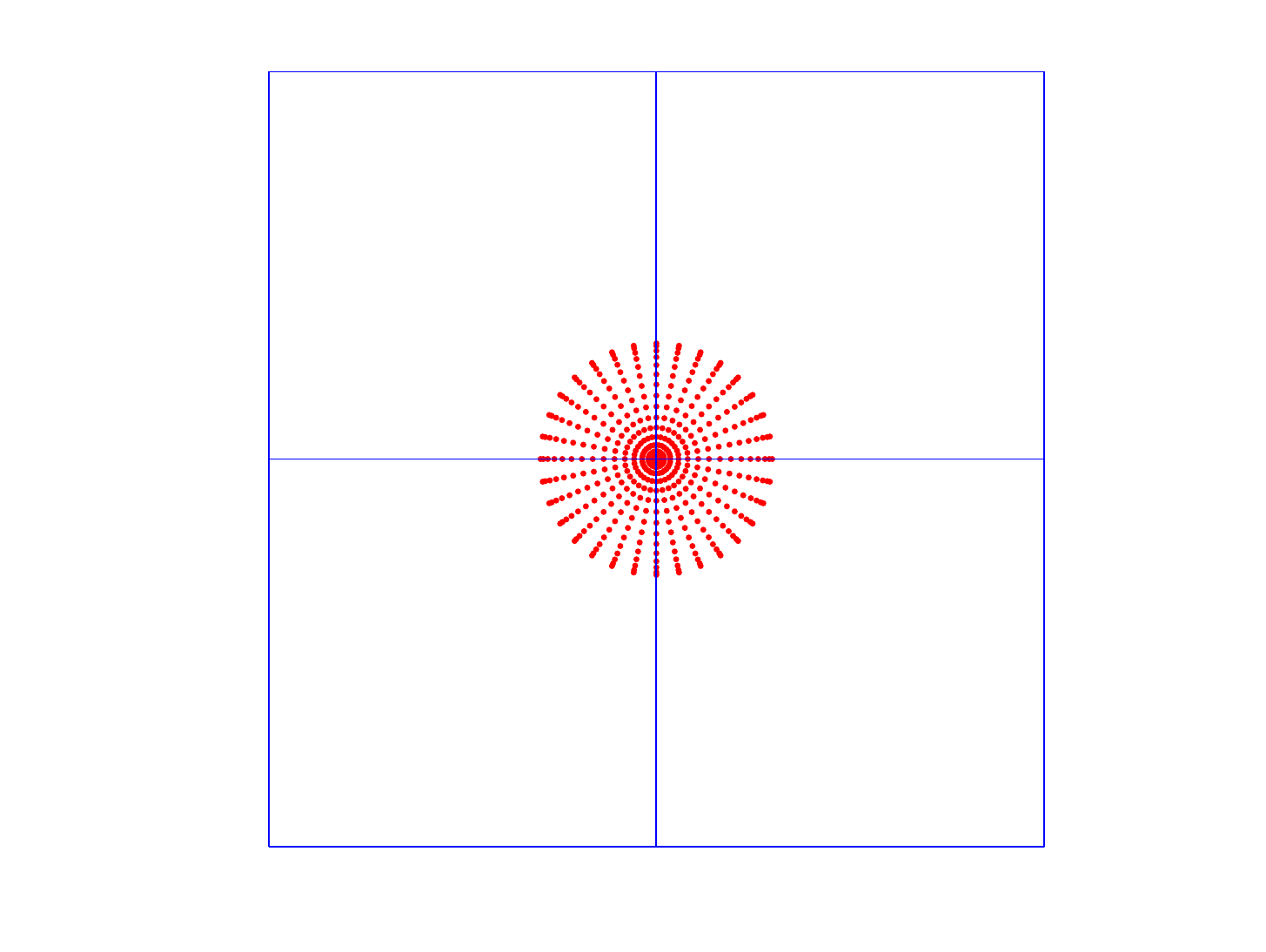}
                \caption{Polar grid}
        \end{subfigure}%
\caption{
Two grids used in the Fourier domain in the improved algorithm in \cite{JGB}: the regular grid on the left panel is used to compute $I_1$ in \eqref{3.9} via the regular FFT;
while the polar grid (confined in a small region centered at the origin)
on the right panel is used to compute
$I_2$ in \eqref{3.9} via the NUFFT. Note that the number of points in the polar grid is
$O(1)$, thus keeping the interpolation cost in NUFFT minimal.
\label{fig3.1}}
\end{figure}

With this $p_d(\bk)$, $I_1$ can be computed via the regular FFT
and $I_2$ can be evaluated via the NUFFT with a fixed (much fewer) number of irregular points
in the Fourier space (see Figure \ref{fig3.1}).
Thus the interpolation cost in the NUFFT is reduced to $O(1)$
and the cost of the overall algorithm is comparable to that of the regular FFT,
with an oversampling factor ($2^3$ for $3$D problems and $2^2$--$3^2$ for $2$D
problems) in front of $O(N\log N)$.

% Remark by Shidong: The following is not correct and commented out by me
% since w_2 is not smooth even though w_2(0)=0!! One can
% only use Gaussian for 3D problems, not 2D problems.

%In many applications, the function $p_d(\bk)$ can be chosen
%as a Gaussian, e.g. $p_d(\bk)=e^{-|\bk|^2/a},\;a>0$ is a positive parameter . In this case,
%$w_3({\bf 0})=\frac{1}{a}$ when $d=3$ and $w_2({\bf 0})=0$ when $d=2$.
%For more details, we
%refer to \cite{JGB} and references therein.

\subsection{Poisson potentials in 2D/1D}
When $U(\bx)$ in (\ref{nonlocal53}) is taken
as the the Green's function of the Laplace operator
$U_{\rm Lap}(\bx)$ (\ref{Kernal21d}) in 2D/1D, the algorithm discussed in the previous section cannot be applied directly to evaluate the Poisson potential
$u(\bx)$  due to the stronger singularity of $\widehat{U}_{\rm Lap}(\bk)=\frac{1}{|\bk|^2}$ at the origin.
Obviously, the Poisson potential $u(\bx)$ satisfies the Poisson equation $-\Delta\, u(\bx) =\rho(\bx)$
with the far field condition
\begin{equation}\label{Poiff35a}
\lim_{|\bx|\to\infty}\left[u(\bx)+\frac{\widehat{\rho}(\bz)}{2\pi} \ln|\bx|\right]=0
\end{equation}
for 2D problems and
\begin{equation}\label{Poiff35b}
\lim_{x\to\pm\infty}\left[u(x)+\frac{1}{2}\left(\widehat{\rho}(0)|x|\mp \widehat{(x \rho)}(0)\right)\right]=0
\end{equation}
for 1D problems, respectively.

Let us first consider the evaluation of the 2D Poisson potential. To overcome the above mentioned difficulties, we  introduce the auxiliary functions
\be\label{GG1}
G(\bx) = \frac{1}{2\pi\sigma^2}e^{- \frac{|\bx|^2}{2\sigma^2}}, \qquad G_1(\bx)=
\widehat{\rho}(\bz) \, G(\bx)  -\, \widehat{(\bx\rho)}(\bz) \cdot \nabla_\bx G(\bx),\qquad
\bx\in {\mathbb R}^2,
\ee
and the function $u_1(\bx)$  which satisfies  the Poisson equation with the far-field condition:
\begin{equation}\label{PoiDiff36}
-\Delta u_1(\bx) =G_1(\bx),\qquad \bx \in \mathbb R^2,
\qquad \quad\lim_{|\bx|\to\infty}\left[u_1(\bx)+\frac{\widehat{\rho}(\bz)}{2\pi} \ln |\bx| \right]=0.
\end{equation}
Here, $\sigma>0$ is a parameter to be chosen later.  Solving (\ref{PoiDiff36}) via the convolution, we
have
\be
u_1(\bx)=(U_{\rm Lap}\ast G_1)(\bx)
=\widehat{\rho}(\bz)\, u_{1,1}(\bx)-\, \widehat{(\bx\rho)}(\bz) \cdot {\bf u}_{1,2}(\bx),
\qquad \bx\in {\mathbb R}^2,
\ee
where
\be\label{u2u3}
u_{1,1}(\bx)=(U_{\rm Lap}\ast G)(\bx), \qquad\qquad {\bf u}_{1,2}(\bx)= \nabla_\bx\, u_{1,1}(\bx),
\qquad \bx\in {\mathbb R}^2.
\ee
Note that $G(\bx)$ is radially symmetric, i.e., $G(\bx)=G(|\bx|)=G(r)$ with $r=|\bx|\ge0$ and
$u_{1,1}(\bx)$  satisfies the Poisson equation
\begin{equation}\label{PoiDiff37}
-\Delta u_{1,1}(\bx) =G(\bx),\qquad \bx \in \mathbb R^2,
\qquad \quad\lim_{|\bx|\to\infty}\left[u_{1,1}(\bx)+\frac{1}{2\pi}\ln |\bx| \right]=0.
\end{equation}
It is clear that $u_{1,1}(\bx)$ is also radially symmetric, i.e., $u_{1,1}(\bx)=u_{1,1}(r)$.
Thus, the Poisson equation (\ref{PoiDiff37}) can be re-formulated as
the following second order ODE:
\be
-\frac{1}{r} \partial_r(r\partial_r u_{1,1}(r))=  G(r), \qquad 0<r<\infty,
\qquad \lim_{r\to\infty} \left[u_{1,1}(r)+ \frac{1}{2\pi}\ln r\right]=0.
\ee
 Integrating the above ODE twice with the far-field boundary condition,
we obtain
\bea\label{phi1_exct}
u_{1,1}(\bx) = \left\{\ba{ll}-\frac{1}{4\pi}
\left[{\textrm E}_1(\frac{|\bx|^2}{2\sigma^2})+2\ln(|\bx|)\right], &\bx\ne {\bf 0},\\[1em]
\frac{1}{4\pi} \left(\gamma_e - \ln(2 \sigma^2) \right), &\bx= {\bf 0},\\
\ea\right. \qquad \bx\in {\mathbb R}^2,
\eea
where ${\textrm E}_1(r):= \int_r^{\infty} t^{-1}e^{-t} {\rm d}t$  for $r>0$
is the exponential integral function \cite{Handbook} and $\gamma_e\approx 0.5772156649015328606$ is the Euler-Mascheroni constant. Differentiating (\ref{phi1_exct}) leads to
\be\label{phi12_ext}
{\bf u}_{1,2}(\bx)= \left\{\ba{ll}
-\frac{1}{2\pi}\frac{\bx}{|\bx|^2}
\left(1-e^{- \frac{|\bx|^2}{2\sigma^2}}\right), &\bx\ne {\bf 0},\\[1em]
0, &\bx= {\bf 0},\\
\ea\right. \qquad \bx\in {\mathbb R}^2.
\ee
Denote
\be \label{uu1u2}
u_2(\bx)=u(\bx)-u_1(\bx)\qquad \Longleftrightarrow \qquad u(\bx)=u_1(\bx)+u_2(\bx), \qquad \bx\in {\mathbb R}^2.
\ee
We have
\begin{equation}
-\Delta u_2(\bx) =\rho(\bx)-G_1(\bx),\qquad \bx \in \mathbb R^2,
\qquad \quad\lim_{|\bx|\to\infty}u_2(\bx)=0.
\end{equation}
Solving the above problem via the Fourier integral, noticing (\ref{GG1}) and
using the fact that
\[\nabla_\bk\widehat{\rho}({\bf 0})=-i\,\widehat{(\bx\rho)}(\bz)=-i\int_{{\mathbb R}^2}\bx
\rho(\bx)\,d\bx,\]
we  obtain
\bea\label{phi2fourier}
u_2(\bx)&=&(U_{\rm Lap}\ast(\rho-G_1))(\bx)
=\frac{1}{(2\pi)^2}\int_{\mathbb R^2} \frac{\widehat{\rho}(\bk)-\widehat{G_1}(\bk)}{|\bk|^2}
\;e^{\,i\; \bk\cdot \bx}\; d \bk\nonumber\\
&=&\frac{1}{(2\pi)^2}\int_{\mathbb R^2} \frac{W(\bk)}{|\bk|}
\;e^{i\, \bk\cdot \bx}\; d \bk\approx
\frac{1}{(2\pi)^2}\int_0^P\int_0^{2\pi}W(\bk)
\;e^{i\, \bk\cdot \bx}\,d|\bk|d\theta, \qquad \bx\in\Omega \subset \mathbb R^2,
\eea
where
\be
W(\bk)=\left\{\ba{ll}\frac{\widehat{\rho}(\bk)-\widehat{G_1}(\bk)}{|\bk|}=
\frac{\widehat{\rho}(\bk)-\big(\widehat{\rho}({\textbf 0})+\bk\cdot \nabla_\bk\widehat{\rho}({\bf 0})
\big)\;e^{-\frac{1}{2} |\bk|^2\sigma^2}}{|\bk|}, &\bk\ne {\bf 0},\\[1em]
0, &\bk={\bf 0},\\
\ea\right.  \qquad \bk \in \mathbb R^2.
\ee
Note that the singularity of $W(\bk)/|\bk|$ at the origin in (\ref{phi2fourier})
is removed by switching to polar coordinates in the Fourier space, and thus
$u_2(\bx)$ can be evaluated by the algorithm in \cite{JGB}.

In practical computations, the parameter $\sigma$ in (\ref{GG1})
should be chosen appropriately such that the Gaussian $e^{-\frac{1}{2} |\bk|^2 \sigma^2}$
and $\bk\cdot \nabla_\bk\widehat{\rho}({\bf 0}) e^{-\frac{1}{2} |\bk|^2 \sigma^2}$
in the Fourier space decay at the same rate or faster than
$\widehat{\rho}(\bk)$ when $|\bk|$ is large.
With this choice of $\sigma$, there is no need to enlarge the computational domain
in the Fourier space for the evaluation of (\ref{phi2fourier}) via the NUFFT.
On the other hand,  there is no need to oversample the
truncated Fourier domain due to the rapid decaying of the Gaussian $e^{-\frac{1}{2} |\bk|^2 \sigma^2}$
in the Fourier space.
Thus, setting the Gaussian to $2\cdot 10^{-16}$ at $|\bk|_\infty=P$ with $P$ being
the side-length of the bounded computational box $B=\{\bk \ |\ |\bk|\le P\} $ in the Fourier space,
we can choose $\sigma =6/P$, a constant that  is independent of the density function $\rho$.

For the convenience of the readers,  we summarize the algorithm to evaluate
the Poisson potential $u(\bx)$  in 2D in Algorithm \ref{alg1}.
\begin{algorithm}
\caption{Evaluation of the Poisson potential in 2D}
\label{alg1}
\begin{algorithmic}[h!]
\STATE Compute $\widehat{\rho}(\bk)$ and $\widehat{(\bx\rho)}({\bf 0})$.
\STATE Evaluate $u_1(\bx) = \widehat{\rho}(\bz)\, u_{1,1}(\bx)-\, \widehat{(\bx\rho)}(\bz) \cdot {\bf u}_{1,2}(\bx)$ via \eqref{phi1_exct} and \eqref{phi12_ext}.
\STATE Evaluate $u_2(\bx)$ through (\ref{phi2fourier}) via the NUFFT \cite{JGB}.
\STATE Compute $u(\bx)=u_1(\bx)+u_2(\bx)$.
\end{algorithmic}
\end{algorithm}

\bigskip
%\vspace{1cm}

Similarly, for the 1D case,  i.e., $U_{\rm Lap}(x)=-\frac{1}{2}|x|$,  we introduce the  auxiliary functions
\bea
G(x) = \frac{1}{\sqrt{2\pi}\,\sigma} e^{-\frac{x^2}{2\sigma^2}},\qquad \qquad
G_1(x) = \widehat{\rho}(0) G(x) - \widehat{(x \rho)}(0)  \,G'(x), \qquad x\in \mathbb R,
\eea
and  function $u_1(x)$ which satisfies the 1D Poisson equation with the far-field condition
\be
\label{1DPois}
-u_1^{\prime\prime}(x) = G_1(x), \qquad x \in {\mathbb R} , \quad \qquad
\lim_{x\to\pm\infty}\left[u_1(x)+ \frac{1}{2}\left(\widehat{\rho}(0)|x|\mp \widehat{(x \rho)}(0)\right) \right]=0.
\ee
Solving the above problem via the convolution, we have
\be
u_1(x)=(U_{\rm Lap}\ast G_1)(x)
=\widehat{\rho}(0)\, u_{1,1}(x)-\, \widehat{(x\rho)}(0)   u_{1,2}(x),
\qquad x\in {\mathbb R},
\ee
where
\bea
&& u_{1,1}(x)=(U_{\rm Lap}\ast G)(x) = -\frac{\sigma}{\sqrt{2\pi} }e^{-\frac{\,x^2}{2 \sigma^2}} -
\frac{1}{2} x \,\textrm{Erf}\left(\frac{\,x\,}{\sqrt{2}\sigma}\right),\\
&&  u_{1,2}(x)= u_{1,1}'(x)= -\frac{1}{2}\,
\textrm{Erf}\left(\frac{\, x\, }{\sqrt{2}\sigma}\right), \qquad x\in {\mathbb R}.
\eea
Here, $ \textrm{Erf}(x) =\frac{2}{\sqrt{\pi}} \int_0^x e^{-t^2} d t $ for $x\in \mathbb R$
is the error function. Combining (\ref{nonlocal53}) and (\ref{1DPois}),
we solve the remaining function $u_2(x) = u(x)-u_1(x)$ via the Fourier integral:
\bea
u_2(x) &=&\left(U_{\rm Lap} \ast (\rho -G_1)\right)(x) = \frac{1}{2\pi} \int_{\mathbb R} \dfrac{\widehat{\rho}(k) - \widehat{G_1}(k)}{k^2} e^{i\,k x} dk \\
&=&\frac{1}{2\pi} \int_{\mathbb R} W(k) e^{i\,k x} dk \approx \frac{1}{2\pi} \int_{-P}^{P} W(k)e^{i\,k x} dk,
\qquad x\in \Omega\subset \mathbb R,
\eea
where
\bea
W(k) = \left\{\ba{ll}\dfrac{\widehat{\rho}(k) - \widehat{G_1}(k)}{k^2} =
\frac{\widehat{\rho}(k)-\big(\widehat{\rho} (0)+k (\widehat{\rho})'
(0)\big)\;e^{-\frac{1}{2} k^2\sigma^2}}{k^2}, & k\ne  0,\\[1em]
-\frac{1}{2} \widehat{(x^2\rho)}(0) + \frac{\sigma^2}{2}\, \widehat{\rho} (0), & k = 0,\\
\ea\right. \qquad k\in \mathbb R.
\eea
Note that the integrand $W(\bk)$ is smooth at the origin $k=0$ in the Fourier space,
therefore $u_2(x)$ can be computed by the regular FFT method. The choice of the parameter $\sigma$
is similar as the one in the 2D case.

We remark that the 1D Poisson potential has also been dealt with successfully in
\cite{SPMCompare} by plugging the Fourier spectral approximation of the density obtained on
a finite interval, e.g., $[-L,L]$, into the convolution \eqref{sps-colb} formula. The method proposed there
is an alternative good choice.

\subsection{Confined Coulomb interactions}
When $U(\bx)$ in (\ref{nonlocal53}) is taken
as the confined Coulomb kernel $U_{\rm Con}^\varepsilon(\bx)$  (\ref{Kernal124d}),
there is no equivalent PDE formulation for the nonlocal potential $u(\bx)$.

When $d=2$, noticing that
\be
\widehat{U}_{\rm Con}^\varepsilon(\bk)\approx \left\{\ba{cl}
\frac{1}{|\bk|}, &|\bk|\to0,\\[0.8em]
\frac{\sqrt{2}}{\sqrt{\pi}\varepsilon |\bk|^2}, &|\bk|\to\infty,\\
\ea\right. \quad \bk \in {\mathbb{R}}^2,
\ee
we can immediately adapt the NUFFT-based solver \cite{JGB} as follows:
\bea\label{ifft2d}
u(\bx)&=&\frac{1}{(2\pi)^{2}}\int_{\mathbb{R}^2}  e^{i\,\bk\cdot\bx}\; \widehat{U}^\varepsilon_{\rm Con}(\bk)\; \widehat{\rho}(\bk)\,d\bk\approx \frac{1}{(2\pi)^{2}}\int_{|\bk|\le P}  e^{i\,\bk\cdot\bx}\; \widehat{U}^\varepsilon_{\rm Con}(\bk)\; \widehat{\rho}(\bk)\,d\bk\nonumber\\
&=&\frac{1}{(2\pi)^{2}}\int_0^P\int_0^{2\pi}
e^{i\,\bk\cdot\bx}\; W_1(\bk) \; \widehat{\rho}(\bk)\,d|\bk|d\theta,
\qquad \bx\in \Omega\subset{\mathbb R}^2,
\eea
where
\be
\label{2Dw1}
W_1(\bk)=|\bk|\,\widehat{U}^\varepsilon_{\rm Con}(\bk)=\frac{2}{\pi}\int_{0}^\infty\frac{|\bk|e^{- \frac{\varepsilon^2s^2}{2}}}
{|\bk|^2+s^2}\,ds=\left\{\ba{ll}
\frac{2}{\pi}\int_{0}^\infty\frac{e^{-\varepsilon^2|\bk|^2s^2/2}}
{1+s^2}\,ds, &\bk\ne {\bf 0}, \\[1em]
1, &\bk= {\bf 0}, \\
\ea\right. \qquad \bk\in {\mathbb R}^2.
\ee
The integral in (\ref{2Dw1}) can be evaluated very accurately via the standard quadrature, such as the Gauss--Kronrod quadrature.

Similarly, when $d=1$ we have
\be
\widehat{U}_{\rm Con}^\varepsilon(k)\approx \left\{\ba{ll}
\frac{1}{2}\left[\ln 2-\gamma_e-2\ln (\varepsilon|k|)\right], &|k|\to0,\\[1em]
\frac{1}{\varepsilon^2 |k|^2}, &|k|\to\infty,\\
\ea\right. \quad k \in {\mathbb{R}}.
\ee
Thus
\bea\label{ifft1d3}
u(x)&=&\frac{1}{2\pi}\int_{\mathbb{R}}
e^{i\,k x}\; \widehat{U}^\varepsilon_{\rm Con}(k)\; \widehat{\rho}(k)\,dk
=-\frac{1}{2\pi}\int_{\mathbb{R}}
e^{i\,k x}k\left[\partial_k\left(\widehat{U}^\varepsilon_{\rm Con}(k)\, \widehat{\rho}(k)\right)
+ix\,\widehat{U}^\varepsilon_{\rm Con}(k)\,\widehat{\rho}(k)\right]\,dk\nonumber\\
&=&-\frac{1}{2\pi}\int_{\mathbb{R}}e^{i\,k x}\left[k\, \partial_k\widehat{U}^\varepsilon_{\rm Con}(k)\, \widehat{\rho}(k)-ik\,\widehat{U}^\varepsilon_{\rm Con}(k)\, \widehat{(x\rho)}(k)
+ixk\, \widehat{U}^\varepsilon_{\rm Con}(k)\, \widehat{\rho}(k)\right]dk\nonumber\\
&=&\frac{1}{2\pi}\int_{\mathbb{R}} e^{i\,k x}\left[W_2(k)\,\widehat{\rho}(k)+i\,
 W_3(k)\,\widehat{(x\rho)}(k)\right]dk
-\frac{i\,x}{2\pi}\int_{\mathbb{R}} e^{i\,k x}\, W_3(k)\,\widehat{\rho}(k)\,dk\nonumber\\
&\approx&\frac{1}{2\pi}\int_{-P}^P e^{i\,k x}\left[W_2(k)\widehat{\rho}(k)+i\,
 W_3(k)\widehat{(x\rho)}(k)\right]dk
-\frac{i\,x}{2\pi}\int_{-P}^P e^{i\,k x}\, W_3(k)\widehat{\rho}(k)\,dk, \
x\in[-L,L].\qquad
\eea
Here
\bea\label{w2w31}
&&W_2(k)=-k\, \partial_k\widehat{U}^\varepsilon_{\rm Con}(k)=
\int_0^\infty
\frac{k^2e^{-\varepsilon^2s/2}}{(k^2+s)^2}\,ds=\left\{\ba{ll}
\int_0^\infty
\frac{e^{-\varepsilon^2k^2s/2}}{(1+s)^2}\,ds, &k\ne 0,\\[1em]
1, &k=0,\\
\ea\right. \quad k\in {\mathbb R}, \\ \label{w2w32}
&&W_3(k)=k\,\widehat{U}^\varepsilon_{\rm Con}(k)=
\int_0^\infty
\frac{k\,e^{-\varepsilon^2s/2}}{2(k^2+s)}\,ds=\left\{\ba{ll}
\int_0^\infty
\frac{k\,e^{-\varepsilon^2k^2s/2}}{2(1+s)}\,ds, &k\ne 0,\\[1em]
0, &k=0,\\
\ea\right. \quad k\in {\mathbb R}.
\eea
The integrals in (\ref{w2w31})-(\ref{w2w32}) can be discretized very accurately via the standard quadrature,
and the integrals in (\ref{ifft1d3}) can be evaluated via the regular FFT.

\medskip

\begin{remark}
If  $\rho(\bx)$ in (\ref{nonlocal53})
is spherically/radially symmetric in 3D/2D, i.e.,
$\rho(\bx)=\rho(|\bx|)=\rho(r)$ with $r=|\bx|$,
and the interaction kernel $U(\bx)$ in (\ref{nonlocal53}) is taken as the Green's function of the Laplace operator
in 3D/2D, then  the nonlocal interaction
$u(\bx)$ in (\ref{nonlocal53}) is also spherically/radially symmetric in 3D/2D,
i.e., $u(\bx)=u(|\bx|)=u(r)$. Additionally, it satisfies the following second-order
ODE
\bea\label{FD-Solver}
&&-\frac{1}{r^{d-1}} \partial_r\left(r^{d-1}\partial_r u(r)\right) =  \rho(r),  \qquad 0 <r <\infty, \qquad d=3,2,\\
\label{FD-Solver2d}
&&\partial_r u(0)=0, \qquad u(r)\to\left\{\ba{ll}0, &d=3,\\[0.5em]
-C_0\,\ln r, &d=2, \\
\ea\right. \qquad r\to\infty,
\eea
where $C_0=\int_0^\infty \rho(r)r\,dr$. Moreover, if $\rho(r)$ has a compact support or
decays exponentially fast when $r\to\infty$, the above problem can be further re-formulated or approximated
by \cite{Han,ExtBd}
\bea\label{FD-Solver2d1}
&&-\frac{1}{r^{d-1}} \partial_r\left(r^{d-1}\partial_r u(r)\right) =  \rho(r),  \qquad 0 <r <L, \qquad d=3,2,\\
\label{FD-Solver2d2}
&&\partial_r u(0)=0, \qquad \partial_r u(L)=\left\{\ba{ll}-\frac{u(L)}{L}, &d=3,\\[1em]
\frac{u(L)}{L\ln L}, &d=2, \\
\ea\right.
\eea
where $L>0$ is large enough such that ${\rm supp}(\rho)\subset [0,L]$ or the truncation
error in $\rho$ outside $[0,L]$ can be negligible. This two-point boundary value problem
can be solved by  the finite difference (FDM) or finite element (FEM) or spectral method.
Comparing to  computing the original convolution or solving the corresponding Poisson equation
in 3D/2D, the memory and/or computational cost are significantly reduced.
\end{remark}

\subsection{Numerical comparisons}
In order to demonstrate the efficiency and accuracy of the NUFFT for the evaluation
of the nonlocal interaction (\ref{nonlocal53}) and compare it with other existing numerical methods,
we adopt the error function
\be
\label{error_infty}
e_h:=\frac{\|u-u_{h}\|_{l^\infty}}{\|u\|_{l^\infty}}=
\frac{\max_{\bx\in \Omega_h} |u(\bx)-u_{h}(\bx)|}{\max_{\bx\in \Omega_h} |u(\bx)|},
\ee
where $\Omega_h$ represents the partition of the bounded computational domain $\Omega$ in 3D/2D
with mesh size $h$, where we usually take  $h_x= h_y=h_z$ in 3D or $h_x= h_y$  in 2D and donate by $h$ unless stated otherwise, and
$u_h(\bx)$ is the numerical solution obtained by a numerical method on the domain $\Omega_h$.
We will compare the method via the NUFFT (referred as {\sl NUFFT}) presented in this section
with those existing numerical methods such as the method via the FFT
(referred as {\sl FFT}) \cite{BMS} and via the DST (referred as {\sl DST})
\cite{DipolarBao,SPMCompare} as well as the finite difference method via (\ref{FD-Solver2d1})-(\ref{FD-Solver2d2})
(referred as {\sl FDM}) \cite{ExtBd} if it is possible.

\medskip
%===============================================================================
%				Example 1
%===============================================================================

{\bf Example 2.1:} {\sl 3D Coulomb interaction.} Here  $d=3$ and $U(\bx)= U_{\rm Cou}(\bx)$,
we take $\rho(\bx):= e^{-(x^{2}+y^2+\gamma^2z^2)/\sigma^2}$ with $\sigma>0$ and $\gamma\ge 1$. The 3D Coulomb interaction can be computed analytically as
\be\label{DensGaus}
u(\bx) =  \left\{\ba{ll}
\frac{\sigma^3 \sqrt{\pi}}{4\;|\bx|\;}\,
\text{Erf}\left( \frac{|\bx|}{ \sigma }\right), &\gamma = 1,\\[1em]
\frac{\sigma^2}{4\gamma} \int_0^\infty \frac{e^{-\frac{x^2+y^2}{\sigma^2(t+1)}}
e^{-\frac{z^2}{\sigma^2(t+\gamma^{-2})}}} {(t+1)\sqrt{t+\gamma^{-2}}}dt     , & \gamma  \ne 1,
\ea\right. \qquad \bx\in{\mathbb R}^3.
\ee

\begin{figure}[t!]
\centerline{ \psfig{figure=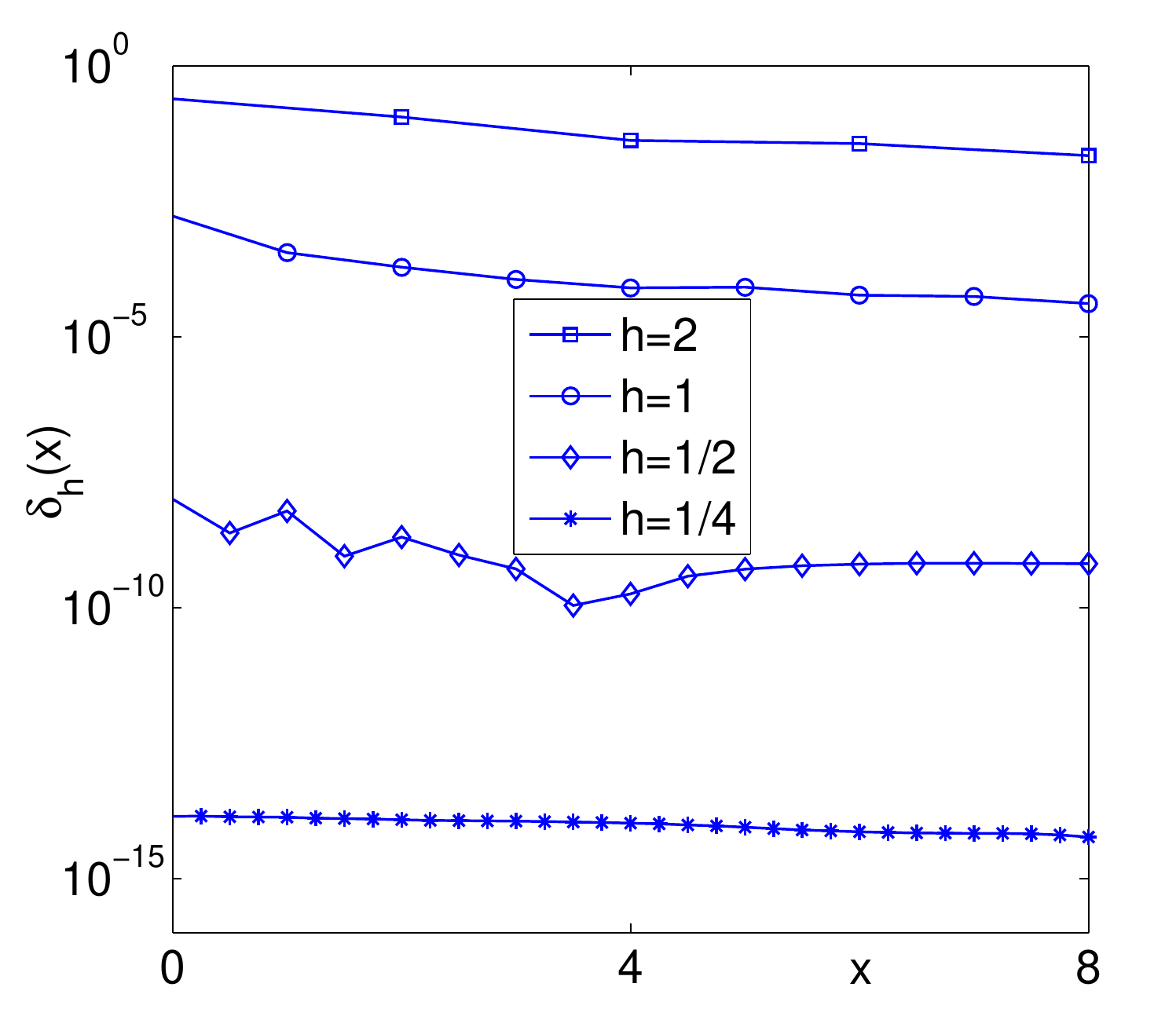,height=6.0cm,width=7.8cm}
\psfig{figure=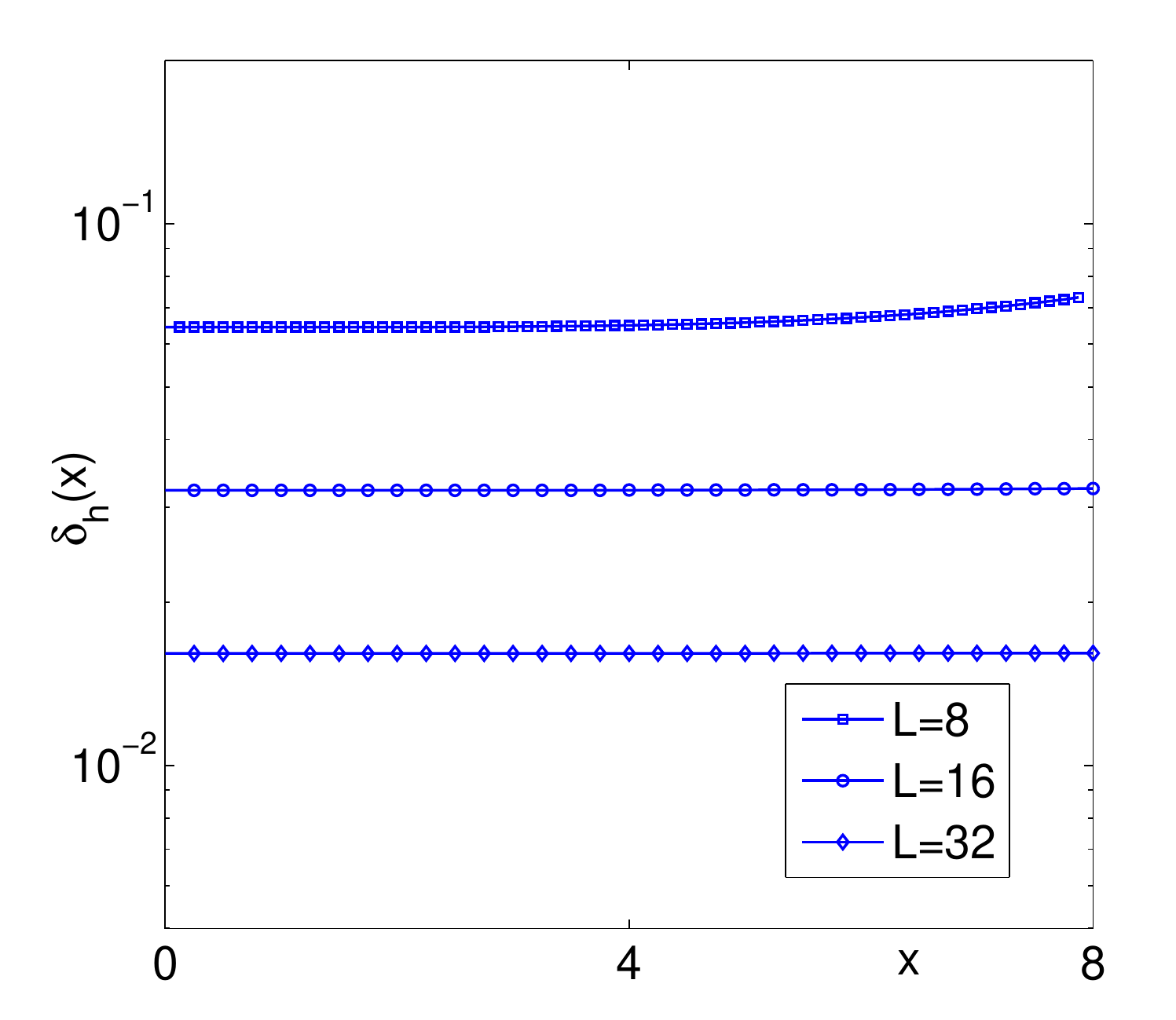,height=6.0cm,width=7.8cm} }
\caption{Errors of  $\delta_h(x)=|u(x,0,0)-u_h(x,0,0)|$ for
the evaluation of the Coulomb interaction in 3D via the NUFFT method
with $L=8$ for different mesh size $h$ (left) and via the DST method
with mesh size $h=1/4$ for different $L$ (right).}
\label{fig:pot:3D:Coulomb}
\end{figure}

The 3D Coulomb interaction $u(\bx)$
is computed numerically via the NUFFT, DST and FFT methods
on a bounded computational domain $\Omega=[-L,L]^2\times [-L/\gamma,L/\gamma]$ with mesh size $h$.
Table \ref{tab:exmp1} shows the errors  $e_h$ via the NUFFT, DST and FFT methods
with $\gamma=1,\sigma = 1.1$ for different mesh size $h$ and  $L$.
Figure \ref{fig:pot:3D:Coulomb} depicts the error of the Coulomb interaction along
the $x$-axis, which is defined as  $\delta_h(x):= |u(x,0,0)-u_h(x,0,0)|$,
obtained via the NUFFT and DST methods with $\gamma=1,\sigma = 1.1$ for different mesh size $h$
and  $L$. In addition, Table  \ref{tab24:exmp1} shows the errors  $e_h$ via
the NUFFT, DST and FFT methods
with $\sigma = 2 $ and $L=8,h=1/4$ for different $\gamma$. Here $h$ denote $h_x = h_y$ and we choose
$h_z = h/\gamma$.

\begin{table}[h!]
\tabcolsep 0pt \caption{Errors for the evaluation of the 3D Coulomb interaction
by different methods for different $h$ and $L$. }
\label{tab:exmp1}
\begin{center}\vspace{-1em}
\def\temptablewidth{1\textwidth}
{\rule{\temptablewidth}{1pt}}
\begin{tabularx}{\temptablewidth}{@{\extracolsep{\fill}}p{1.35cm}rllll}
NUFFT & $h = 2$& $h = 1$  & $h= 1/2$ & $h=1/4$ & $h  = 1/8$  \\[0.25em]\hline
$L=4$ &4.191E-01 & 2.696E-03 & 6.634E-07 & 4.599E-07&  3.688E-07\\
 $L = 8$ &4.111E-01&2.817E-03 & 1.667E-08 & 2.367E-14&2.404E-14\\
 $L = 16$&4.127E-01&2.848E-03& 1.732E-08 & 1.420E-14 &1.334E-14 \\
\hline\hline
DST & $h = 2$& $h = 1$  & $h= 1/2$ & $h=1/4$ & $h  = 1/8$  \\[0.25em]\hline
$L=4 $   &2.437E-01&2.437E-01 &2.437E-01 &2.437E-01 &2.437E-01\\
$L= 8$  &2.754E-01&1.219E-01 &1.219E-01 &1.219E-01 &1.219E-01\\
$L=16$  &3.433E-01&6.093E-02 &6.093E-02 &6.093E-02 &6.093E-02\\
$L=32$  &3.780E-01&3.046E-02 &3.046E-02 &3.046E-02 &3.046E-02\\
$L=64$  &3.956E-01&1.523E-02 &1.523E-02 &1.523E-02 &1.523E-02\\
\hline\hline
FFT & $h = 2$& $h = 1$  & $h= 1/2$ & $h=1/4$ & $h  = 1/8$  \\[0.25em]\hline
$L=4$  &3.032E-01 &3.363E-01 &3.385E-01 &3.385E-01 &3.385E-01\\
$L= 8$ &1.744E-01 &1.712E-01 &1.720E-01 &1.720E-01 &1.720E-01\\
$L=16$ &2.958E-01 & 8.666E-02 & 8.632E-02 & 8.632E-02 & 8.632E-02\\
$L=32$ &3.550E-01 & 4.372E-02 & 4.320E-02 & 4.320E-02 & 4.320E-02\\
$L=64$ &3.843E-01 & 2.214E-02 & 2.161E-02 & 2.161E-02 & 2.161E-02\\
\end{tabularx}
{\rule{\temptablewidth}{1pt}}
\end{center}
\end{table}

\begin{table}[h!]
\tabcolsep 0pt \caption{Errors for the evaluation of the 3D  Coulomb interaction
by different methods with $\sigma = 2$ and $L=8,h=1/4$ for different $\gamma$. }
\label{tab24:exmp1}
\begin{center}\vspace{-1em}
\def\temptablewidth{1\textwidth}
{\rule{\temptablewidth}{1pt}}
\begin{tabularx}{\temptablewidth}{@{\extracolsep{\fill}}p{1.35cm}rlll}
 & $\gamma = 1$& $\gamma = 2$  & $\gamma= 4$ & $\gamma=8$   \\[0.25em]\hline
NUFFT &2.164E-14 &2.134E-14 &2.044E-14  &2.005E-14 \\ \hline
DST &0.146 &0.441 &1.559 &3.782\\ \hline
FFT  &0.208 &0.310 &1.327 &3.349
\end{tabularx}
{\rule{\temptablewidth}{1pt}}
\end{center}
\end{table}

From Tables \ref{tab:exmp1}--\ref{tab24:exmp1} and Figure \ref{fig:pot:3D:Coulomb}, we can observe clearly that :
(i)  The errors are saturated in the DST and FFT methods as mesh size $h$
tends smaller and the saturated accuracies decrease linearly with respect to the box size $L$;
(ii) The NUFFT method is spectrally accurate and it essentially does not depend on the domain,
which implies that a very large bounded computational domain is not necessary in practical computations
when the NUFFT method is used;
(iii) The NUFFT is capable of dealing with anisotropic densities, which is quite useful in numerical
simulation of BEC with strong confinement,
while the errors by the DST and FFT methods increase dramatically with strongly anisotropic densities
 (cf. Tab. \ref{tab24:exmp1}).

\medskip

%===============================================================================
%				Example 2
%===============================================================================

{\bf Example 2.2:} {\sl 2D Coulomb interaction.} Here $d=2$ and $U(\bx)= U_{\rm Cou}(\bx)$,
we take  $\rho(\bx):= e^{-(x^{2}+\gamma^2y^2)/\sigma^2}$ with $\sigma>0$ and $\gamma\ge1$ .
The 2D Coulomb interaction can be obtained analytically as
\be\label{2.5D-clb-exact}
u(\bx) =  \left\{\ba{ll} \frac{\sqrt{\pi}\, \sigma}{2}
\,{\mathrm I}_0\left(\frac{|\bx|^2}{2 \sigma^2}\right)\,e^{-\frac{|\bx|^2}{2\sigma^2}}, &\gamma = 1,\\[0.45em]
\frac{\sigma}{\gamma \sqrt{\pi}} \int_0^\infty \frac{e^{-\frac{x^2}{\sigma^2(t^2+1)}}
e^{-\frac{y^2}{\sigma^2(t^2+\gamma^{-2})}}} {\sqrt{t^2+1}\sqrt{t^2+\gamma^{-2}}}dt     , & \gamma  \ne 1,
\ea\right.\qquad \bx\in{\mathbb R}^2,
\ee
where ${\mathrm I}_0$ is the modified Bessel function of order zero \cite{Handbook}.
To numerically compute the integral  in \eqref{2.5D-clb-exact}, we first split it into two integrals
and reformulate the one with infinite interval into some equivalent integral with finite interval by a simple change of variable.
We then apply the Gauss--Kronrod quadrature to each with fine accuracy control so as to achieve accurate reference solutions.

The 2D Coulomb interaction $u(\bx)$ is computed numerically via the NUFFT, DST and FFT methods
on a bounded computational domain $\Omega=[-L,L]\times[-L/\gamma,L/\gamma]$ with mesh size $h$.
Table  \ref{tab:sqrt:nufft} shows the errors  $e_h$ via the NUFFT, DST and FFT methods with $\sigma= \sqrt{1.2}$ and $\gamma=1$ under different mesh size $h$ and  $L$. In addition, Table  \ref{tab34:exmp1} shows the errors  $e_h$ via the NUFFT, DST and FFT methods with $\sigma = 2$,  $L=12$ and $h=1/8$ for different $\gamma$.Here $h$ denote $h_x$ and we choose
$h_y = h/\gamma$.

\begin{table}[h!]
\tabcolsep 0pt \caption{Errors for the evaluation of the 2D Coulomb interaction
by different methods for different $h$ and $L$.}
\label{tab:sqrt:nufft}
\begin{center}\vspace{-0.5em}
\def\temptablewidth{1\textwidth}
{\rule{\temptablewidth}{1pt}}
\begin{tabularx}{\temptablewidth}{@{\extracolsep{\fill}}p{1.25cm}llllll}
NUFFT& $h = 2$ & $h = 1$  & $h =1/2$ & $h=1/4$ & $h = 1/8$\\ \hline
$L=4$  & 1.837 &  5.540E-02  & 4.289E-07  & 3.383E-07 & 2.937E-07\\
$L=8$ &4.457E-01 &2.373E-03 &2.714E-08 &3.202E-15 &3.431E-15\\
$L=16 $ &2.084E-01 &2.385E-03 &2.761E-08 &2.745E-15  &2.859E-15\\ \hline \hline
DST &$h = 2$ & $h = 1$  & $h =1/2$ & $h=1/4$ & $h = 1/8$ \\  \hline
$L=4$   &1.577E-01 &1.577E-01 &1.577E-01 &1.577E-01 &1.577E-01\\
$L = 8$ &1.348E-01 &7.762E-02 &7.762E-02 &7.762E-02 &7.762E-02\\
$L = 16$&1.711E-01 &3.867E-02 &3.867E-02 &3.867E-02 &3.867E-02\\
$L =32$ &1.897E-01 &1.932E-02 &1.932E-02 &1.932E-02 &1.932E-02\\
$L=64$  &1.991E-01 &9.658E-03 &9.658E-03 &9.658E-03 &9.658E-03 \\
\hline \hline
FFT &$h = 2$ & $h = 1$  & $h =1/2$ & $h=1/4$ & $h = 1/8$ \\  \hline
$L=4$   &2.855E-01 &2.961E-01 &2.980E-01 &2.980E-01 &2.980E-01\\
$L = 8$ &1.553E-01 &1.503E-01 &1.502E-01 &1.502E-01 &1.502E-01\\
$L = 16$&1.157E-01 &7.596E-02 &7.528E-02 &7.528E-02 &7.528E-02\\
$L =32$ &1.624E-01 &3.843E-02 &3.766E-02 &3.766E-02 &3.766E-02\\
$L=64$  &1.856E-01 &1.961E-02 &1.883E-02 &1.883E-02 &1.883E-02\\
\end{tabularx} {\rule{\temptablewidth}{1pt}}
\end{center}
\end{table}

\begin{table}[h!]
\tabcolsep 0pt \caption{Errors for the evaluation of the 2D  Coulomb interaction
by different methods with $L=12,h=1/8$ for different $\gamma$. }
\label{tab34:exmp1}
\begin{center}\vspace{-1em}
\def\temptablewidth{1\textwidth}
{\rule{\temptablewidth}{1pt}}
\begin{tabularx}{\temptablewidth}{@{\extracolsep{\fill}}p{1.35cm}rlll}
 & $\gamma = 1$& $\gamma = 2$  & $\gamma= 4$ & $\gamma=8$   \\[0.25em]\hline %$\gamma  = 16$
 NUFFT &4.230E-14 &3.102E-15 &3.504E-15 &4.381E-15 \\ \hline %&8.073E-15
DST   &0.373 &0.386 &0.412 &0.446\\ \hline% &0.488
FFT   &0.426 &0.425 &0.405 &0.344% &0.231
\end{tabularx}
{\rule{\temptablewidth}{1pt}}
\end{center}
\end{table}

From Tables \ref{tab:sqrt:nufft}-\ref{tab34:exmp1}, we can conclude that:
(i) The errors obtained by the DST and FFT methods reach a saturation accuracy on any fixed
domain and  we can observe a first order convergence in the saturated accuracy with
respect to the domain size  $L$.
(ii)  The NUFFT method is spectrally  accurate  and it essentially does not  depend
on the domain which makes it perfect for computing the whole space potential.
(iii) The NUFFT is capable of dealing with anisotropic densities,
while the results obtained by the  DST and FFT methods are far from the exact solutions
when the bounded computational domain is not large enough.

\medskip

%===============================================================================
%				Example 1
%===============================================================================

{\bf Example 2.3:} {\sl 2D Poisson potential.} Here $d=2$ and $U(\bx)= U_{\rm Lap}(\bx)$,
we take $\rho(\bx):= e^{-|\bx|^2/\sigma^2}=e^{-r^2/\sigma^2}$ with $r=|\bx|$ and
$\sigma>0$. The 2D Poisson potential can be obtained analytically as
\be\label{eg-2d-Exact}
u(\bx)=-\frac{\sigma^2}{4}\,\left[{\textrm E}_1\left( \frac{|\bx|^2}{\sigma^2}\right)+2\ln(|\bx|)\right],
\qquad \bx\in{\mathbb R}^2.
\ee

In this case, we choose $\sigma=\sqrt{1.3}$. The 2D Poisson potential $u(\bx)$
is computed numerically via the NUFFT method
on a bounded computational domain $\Omega=[-L,L]^2$ with mesh size $h$
and the FDM through the formulation (\ref{FD-Solver2d1})-(\ref{FD-Solver2d2})
on the interval $[0,L]$ with mesh size $h$.

Table \ref{tab_NUFFT_2D_Poisson} shows the errors of the 2D Poisson potential obtained by the
NUFFT solver on a square domain and the errors by the FDM solver as well as its
convergence rate with respect to the mesh size $h$. In addition, to demonstrate the efficiency of the NUFFT method, Table \ref{tab:2d:Poisson:time}  displays the computational time (CPU time in seconds) of the NUFFT solver with $L=16$ and $h=1/4$, where the time is measured when the algorithm is implemented in Fortran, the code is compiled by ifort 13.1.2 using the option -g, and executed on 32-bit
Ubuntu Linux on a 2.90GHz Intel(R) Core(TM) i7-3520M CPU with 6MB cache.

\begin{table}[h!]
\tabcolsep 0pt \caption{Errors for the evaluation of the 2D Poisson potential
 by different methods for different $h$ and $L$.}
\label{tab_NUFFT_2D_Poisson}
\begin{center}\vspace{-0.5em}
\def\temptablewidth{1\textwidth}
{\rule{\temptablewidth}{1pt}}
\begin{tabularx}{\temptablewidth}{@{\extracolsep{\fill}}p{1.25cm}rlllll}
NUFFT &  $h = 2$  & $h =1$ & $h=1/2$ & $h = 1/4$ & $h=1/8$ \\[0.25em]\hline
$L=4$ &5.821E-01 &1.133E-02 &3.011E-06 &1.994E-06 &1.650E-06\\
$L=8$ &1.685E-01 & 6.820E-04 &1.754E-09 &4.936E-14 &4.857E-14\\
$L=16$&1.684E-01 & 5.333E-04 &1.391E-09 &4.577E-14 &4.561E-14\\ \hline\hline
FDM &$h = 1/4$  & $h = 1/8$ & $h=1/16$ & $h = 1/32$ & $h = 1/64$\\[0.25em]\hline% &$h = 1/128$
$L=4$         &4.646E-03 & 1.155E-03 &  2.910E-04 & 7.602E-05 &  2.246E-05 \\ % & 9.094E-06\\
{\textrm rate}&- & 2.0081 &    1.9889  &   1.9365   &  1.7590  \\ % &  1.3043\\
 $L=8$         &4.101E-03  & 1.019E-03   & 2.542E-04   & 6.353E-05   & 1.588E-05 \\ %&    3.970E-06\\
 {\textrm rate}&- &2.0093 &   2.0024&    2.0006   & 2.0002\\ %  & 2.0000 \\
  $L=16$       &4.052E-03  &  1.007E-03   & 2.512E-04  &  6.278E-05   & 1.569E-05 \\ % &  3.923E-06\\
 {\textrm rate}&- &2.0092 &   2.0023   & 2.0006  &  2.0001 %  & 2.0000
\end{tabularx}
{\rule{\temptablewidth}{1pt}}
\end{center}
\end{table}

\begin{table}[h!]
\tabcolsep 0pt
\begin{center}
\caption{CPU time (in seconds) of  the NUFFT solver for the evaluation of the 2D Poisson potential.
Here $T_{\rm FFT}$  and  $T_{\rm NUFFT}$ are the time for
the evaluation of $I_1$ and $I_2$ in (\ref{3.9}) via the FFT and NUFFT methods,
respectively. }
\label{tab:2d:Poisson:time}
\vspace{-0.5em}
\def\temptablewidth{1\textwidth}
{\rule{\temptablewidth}{1pt}}
\begin{tabularx}{\temptablewidth}{@{\extracolsep{\fill}}p{1.25cm}ccc}
 & $T_{\rm FFT}$ & $T_{\rm NUFFT}$ & $T_{\rm Total}$ \\ \hline
$h  = 1$  & 0.01 &0.05 &0.06 \\
$h  = 1/2$& 0.02 &0.08 &0.10\\
$h  = 1/4$& 0.12 &0.20 &0.32\\
$h\!=1/8\!$&0.60 &0.78 &1.38
\end{tabularx}
{\rule{\temptablewidth}{1pt}}
\end{center}
\end{table}

From Tables \ref{tab_NUFFT_2D_Poisson}--\ref{tab:2d:Poisson:time}, we can see clearly that:
(i) The NUFFT solver is spectrally accurate while the FDM solver is only second order accurate,
and the NUFFT solver is much more accurate than the FDM solver.
(ii) The errors obtained by both methods do not essentially depend on the domain size;
(iii) The complexity of the NUFFT solver scales like $O(N\ln N)$ as expected, which is the same as  those presented in \cite{JGB}.

\medskip

%==================================================================
%				new section
%==================================================================

\section{Computing the ground state}\setcounter{equation}{0}
In this section, we present an efficient and accurate
numerical method for computing the ground state of
(\ref{ground}) by combining  NUFFT-based nonlocal interaction potential solver
and the normalized gradient flow that is discretised by backward Euler Fourier pseudospectral
method, and compare it with those existing numerical methods.

\subsection{A numerical method via the  NUFFT}
We choose $\tau>0$ as the time step and denote $t_n=n\tau$ for $n=0,1,2,\ldots$\,.
Different efficient and accurate numerical methods have
been proposed in the literature  for computing
the ground state \cite{DipJCP,BCL_JCP,BD_SISC,Dong,SPMCompare}. One of the most simple and popular
methods is through the following gradient flow with discretized normalization
(GFDN):
\bea \label{gf-eq1}
&&\partial_t \phi(\bx,t)= \left[\frac{1}{2} \Delta  -V(\bx)- \beta \,\varphi(\bx,t) \right]
\phi(\bx,t),\qquad \bx \in {\mathbb{ R}}^d,\quad t_n\le t<t_{n+1},\\
\label{gf-colb1}
&&\quad \varphi(\bx,t) = \left(U\ast |\phi|^2\right)(\bx,t), \qquad \qquad
\bx \in {\mathbb{ R}}^d, \quad t_n\le t< t_{n+1},\\
\label{gffg3}
&&\phi(\bx,t_{n+1}):=\phi(\bx,t_{n+1}^+)=\frac{\phi(\bx,t_{n+1}^-)}{\|\phi(\bx,t_{n+1}^-)\|}, \qquad
\bx \in {\mathbb{ R}}^d, \qquad n=0,1,2,\ldots\;
\eea
with the initial data
\be
\phi(\bx,0)=\phi_0(\bx), \qquad \bx \in {\mathbb{ R}}^d, \qquad {\rm with} \qquad \|\phi_0\|^2:=\int_{{\mathbb R}^d}
|\phi_0(\bx)|^2\,d\bx=1.
\ee
Let $\phi^n(\bx)$ and $\varphi^n(\bx)$ be the numerical approximation of $\phi(\bx,t_n)$ and
$\varphi(\bx,t_n)$, respectively, for $n\ge0$.
The above GFDN is usually discretized in time via the backward Euler method \cite{DipJCP,BCL_JCP,BD_SISC,Dong,SPMCompare}
\bea \label{gf-eq2}
&&\frac{\phi^{(1)}(\bx)-\phi^{n}(\bx)}{\tau}= \left[\frac{1}{2} \Delta
-V(\bx)- \beta \,\varphi^n(\bx) \right]
\phi^{(1)}(\bx),\qquad \bx \in {\mathbb{ R}}^d,\\
\label{gf-colb2}
&&\quad \varphi^n(\bx) = \left(U\ast |\phi^n|^2\right)(\bx), \qquad \qquad
\bx \in {\mathbb{ R}}^d, \\
\label{gffg33}
&&\phi^{n+1}(\bx)=\frac{\phi^{(1)}(\bx)}{\|\phi^{(1)}(\bx)\|}, \qquad
\bx \in {\mathbb{ R}}^d, \qquad n=0,1,2,\ldots\;.
\eea
Then an efficient and accurate numerical method can be designed
by: (i) truncating the above problem on a bounded computational domain $\Omega$ with
periodic BC on $\partial\Omega$;
(ii) discretizing in space via the Fourier pseudospectral method;
and (iii) evaluating the nonlocal interaction
$\varphi^n(\bx)$ in (\ref{gf-colb2})
by the algorithm via the NUFFT discussed in the previous section.
When $\phi_0(\bx)$ is chosen as a positive function, the ground state
can be obtained as $\phi_g(\bx)=\lim_{n\to\infty}\phi^n(\bx)$ for $\bx\in \Omega$.
The details are omitted here for brevity and this method is referred as
the {\sl GF-NUFFT} method. We remark here that $\widehat{|\phi^n|^2}({\bf 0})=1$
for $n\ge0$.

For comparison, for the Coulomb interaction in 3D/2D,
when the NUFFT solver is replaced by the standard FFT,
we refer the method as {\sl GF-FFT}. In addition,
when  (\ref{gf-colb2}) is reformulated as its equivalent PDE formulation
(\ref{PoiDiff})-(\ref{SqrtPoiDiff})  on  $\Omega$
with homogeneous Dirichlet  BC on $\partial\Omega$ and
 solved via the sine pseudospectral method
 \cite{DipJCP,BJNY,SPMCompare}, we refer it as {\sl GF-DST}.

\subsection{Numerical comparisons}

In order to compare the {\sl GF-NUFFT} method with   {\sl GF-FFT} and {\sl GF-DST} methods for
computing the ground state, we denote $\varphi_g(\bx)=(U\ast\phi_g)(\bx)$ and introduce
the errors
\[e_{\phi_g}^h:=\frac{\max_{\bx\in \Omega^h}|\phi_g(\bx)-\phi_g^h(\bx)|}{\max_{\bx\in \Omega^h}|\phi_g(\bx)|} ,
\qquad e_{\varphi_g}^h:=\frac{\max_{\bx\in \Omega^h} |\varphi_g(\bx)-\varphi_g^h(\bx)|}{\max_{\bx\in \Omega^h}|\varphi_g(\bx)|},
\]
where $\phi_g^h$ and $\varphi_g^h$ are obtained numerically by a numerical method
with mesh size $h$. Additionally, we split the energy functional into three parts
\[E(\phi)=E_{\rm kin}(\phi)+E_{\rm pot}(\phi)+E_{\rm int}(\phi),\]
where the   kinetic energy  $E_{\rm kin}(\phi)$, the potential energy  $E_{\rm pot}(\phi)$ and the interaction
energy $E_{\rm int}(\phi)$ are defined as
\[E_{\rm kin}(\phi)   =  \frac{1}{2} \int_{\mathbb R^d} |\nabla \phi(\bx)|^2 d \bx,\qquad
E_{\rm pot}(\phi)  =  \int_{\mathbb R^d} V(\bx) |\phi(\bx)|^2 d \bx, \qquad
E_{\rm int}(\phi)  =  \frac{\beta}{2}\int_{\mathbb R^d}\varphi(\bx) |\phi(\bx)|^2  d\bx,
\]
respectively. Moreover, the chemical potential can be reformulated as $\mu(\phi)=E(\phi)+E_{\rm int}(\phi)$. Furthermore, if   the external potential $V(\bx)$ in (\ref{sps-eq}) was taken as the harmonic potential \cite{Bao2013,BJNY,ExtBd}, the energies of the ground state satisfy  the following viral identity
\[
0=I:=2 E_{\rm kin}(\phi_g)  - 2 E_{\rm pot}(\phi_g)  +\left\{\ba{ll}
E_{\rm int}(\phi_g), &U=U_{\rm Cou} \ \hbox{in 3D/2D},\\[0.5em]
\frac{\beta}{4\pi}, &U=U_{\rm Lap}\ \hbox{in 2D}.\\
\ea\right.
\]
We denote $I^h$ as an approximation of $I$ when $\phi_g$ is replace by $\phi_g^h$ in the above equality.
In our computations, the ground state $\phi_g^h$ is reached numerically when $\max_{\bx\in \Omega^h}\frac{|\phi^{n+1}(\bx)-\phi^n(\bx)|}{\tau}\le \varepsilon_0$ with $\varepsilon_0$  a prescribed accuracy, e.g., $\varepsilon_0=10^{-10}$.  The initial data $\phi_0(\bx)$ is chosen as a Gaussian and the time step is taken as $\tau=10^{-2}$.  In the comparisons, the  ``exact" solution $\phi_g(\bx)$ was
obtained numerically via the GF-NUFFT method on a  large enough domain $\Omega$ with small enough mesh size $h$ and time step $\tau$.

%========================================================================
%			GS-Example 3
%========================================================================

\medskip

{\bf Example 3.1:} {\sl The NLSE with the Coulomb interaction in 3D.} We take $d=3$ and
$U(\bx)= U_{\rm Cou}(\bx)$
in (\ref{sps-eq})-(\ref{sps-colb}). The ground state is computed numerically on
a bounded domain $\Omega=[-8,8]^3$.  Table \ref{tab:gs:SPS:3d} shows
the errors $e_{\phi_g}^h$ and $e_{\varphi_g}^h$ with $V(\bx) = \frac{1}{2}(x^2+y^2+z^2)$
in (\ref{sps-eq}) for
different numerical methods, $\beta$ and mesh size $h$.
In addition, Table \ref{tab:virial:3d:sps}
lists the energy $E_g:=E(\phi_g^h)$, chemical potential $\mu_g:=\mu(\phi_g^h)$,
kinetic energy $E_{\rm kin}^g:=E_{\rm kin}(\phi_g^h)$, potential energy $E_{\rm pot}^g:=E_{\rm pot}(\phi_g^h)$,
interaction energy $E_{\rm int}^g:=E_{\rm int}(\phi_g^h)$
and  $I^h$ with $h=1/8$ and $V(\bx) = \frac{1}{2}(x^2+y^2+4z^2)$
in (\ref{sps-eq}) for different $\beta$.

%%%%%%% Ground state of the 3D SPS by Backward Euler + NUFFT solver method %%%%%%%%
\begin{table}[h!]
\tabcolsep 0pt \caption{Errors of the ground state for the NLSE with the
3D Coulomb interaction for different methods and mesh size $h$. }
\label{tab:gs:SPS:3d}
\begin{center}\vspace{-1em}
\def\temptablewidth{1\textwidth}
{\rule{\temptablewidth}{1pt}}
\begin{tabularx}{\temptablewidth}{@{\extracolsep{\fill}}p{0.025cm}lllll}
\multicolumn{2}{c|}{GF-NUFFT} & $h=2$ &  $h = 1$  & $h =1/2$ & $h=1/4$  \\ \hline
\multirow{2}*{\scalebox{1.25}{$e^h_{\phi_g}$}}&\multicolumn{1}{|c|}{$\beta=-5\;$}  & 5.362E-02 &  1.954E-04 &  2.201E-07 & 4.643E-11  \\
&\multicolumn{1}{|c|}{$\beta=\;\;5\;$}& 1.512E-01 &  4.712E-04 &  4.026E-08 &  1.141E-10  \\ \hline
\multirow{2}*{\scalebox{1.25}{$e^h_{\varphi_g}$}} &\multicolumn{1}{|c|}{$\beta=-5\;$}   & 2.532E-01 &  3.769E-03 &  8.153E-07 &  7.035E-11  \\
&\multicolumn{1}{|c|}{$\beta=\;\;5\;$}& 2.682E-01 &  7.061E-04 &  1.225E-07 &  8.048E-11  \\
\hline\hline
\multicolumn{2}{c|}{GF-DST}&$h=2$ &  $h = 1$  & $h =1/2$ & $h=1/4$ \\ \hline %& $h =1/8$
\multirow{2}*{\scalebox{1.25}{$e^h_{\phi_g}$}}&\multicolumn{1}{|c|}{$\beta=-5\;$} &2.319E-01 &  9.439E-03 &  1.637E-06 &  6.309E-07 \\ %&  6.309E-07
&\multicolumn{1}{|c|}{$\beta=\;\;5\;$}&1.659E-01 &  9.469E-04 &  8.306E-07 &  8.531E-07 \\ \hline % &  8.531E-07
\multirow{2}*{\scalebox{1.25}{$e^h_{\varphi_g}$}}&\multicolumn{1}{|c|}{$\beta=-5\;$}&7.297E-02 &  9.551E-02 &  9.945E-02 &  1.027E-01  \\ %&  1.049E-01
&\multicolumn{1}{|c|}{$\beta=\;\;5\;$}&7.809E-02 &  1.016E-01 &  1.057E-01 &  1.091E-01  \\ %&  1.114E-01
\end{tabularx}
{\rule{\temptablewidth}{1pt}}
\end{center}
\end{table}

\begin{table}[htbp]
\tabcolsep 0pt \caption{Different energies of the ground state and $I^h$ for the NLSE with the
3D Coulomb interaction  for different $\beta$.}
\label{tab:virial:3d:sps}
\begin{center}\vspace{-0.5em}
\def\temptablewidth{1\textwidth}
{\rule{\temptablewidth}{1pt}}
\begin{tabularx}{\temptablewidth}{@{\extracolsep{\fill}}p{1.25cm}lccccc}
$\beta$  &$E_g$ &$\mu_g$ &$E_{\textrm {kin}}^g$  & $E_{\textrm {pot}}^g$ & $E_{\textrm {int}}^g$ &$ I^h $ \\
\hline
$-10$ & 1.6370& 1.2630& 1.0990  &  9.1197E-01 & -3.7401E-01& -3.39E-10\\
$-5$  & 1.8212& 1.6397& 1.0467  &  9.5594E-01 & -1.8147E-01& -3.63E-10\\
$-1$  & 1.9646& 1.9292& 1.0089  &  9.9118E-01 & -3.5462E-02& -3.87E-10\\
$ 1$  & 2.0351& 2.0702& 9.9128E-01 &  1.0088  &  3.5064E-02& -3.86E-10\\
$5$   & 2.1739& 2.3454& 9.5831E-01 &  1.0441  &  1.7151E-01& -4.30E-10\\
$10$  & 2.3431& 2.6772& 9.2101E-01 &  1.0880  &  3.3408E-01& -1.16E-10\\
\end{tabularx}
{\rule{\temptablewidth}{1pt}}
\end{center}
\end{table}

\
%========================================================================
%			GS-Example 2
%========================================================================

\
{\bf Example 3.2:} {\sl The NLSE with the Coulomb interaction in 2D.} We take $d=2$ and
$U(\bx)= U_{\rm Cou}(\bx)$ in (\ref{sps-eq})-(\ref{sps-colb}).
The ground state is computed numerically on
a bounded domain $\Omega=[-L,L]^2$ with different  mesh size $h$.  Table \ref{gs_2d_Coulomb} shows
the errors $e_{\phi_g}^h$ and $e_{\varphi_g}^h$ with $V(\bx) = \frac{1}{2}(x^2+4y^2) $  for different numerical methods, $\beta$ and mesh size $h$ on $[-L,L]^2$.
In addition, Table \ref{virial_2d_Coulomb}
lists the energy $E_g:=E(\phi_g^h)$, chemical potential $\mu_g:=\mu(\phi_g^h)$,
kinetic energy $E_{\rm kin}^g:=E_{\rm kin}(\phi_g^h)$, potential energy $E_{\rm pot}^g:=E_{\rm pot}(\phi_g^h)$,
interaction energy $E_{\rm int}^g:=E_{\rm int}(\phi_g^h)$
and $I^h$ with $h=1/8$ and $V(\bx) =\frac{1}{2}(x^2+4y^2)$  on $[-8,8]^2$ for different $\beta$.  \\

%%%%%%%    NLS with 2D Coulomb interaction  %%%%%%%%
\begin{table}[h!]
\tabcolsep 0pt \caption{Errors of the ground state for the NLSE with 2D
Coulomb interaction on $[-L,L]^2$ with mesh size $h$. } \label{gs_2d_Coulomb}
\begin{center}\vspace{-1em}
\def\temptablewidth{1\textwidth}
{\rule{\temptablewidth}{1pt}}
\begin{tabularx}{\temptablewidth}{@{\extracolsep{\fill}}p{0.025cm}lllll}
\multicolumn{2}{c|}{GF-NUFFT ($L =8$)}  &  $h = 1$  & $h =1/2$ & $h=1/4$ & $h = 1/8$ \\ \hline
\multirow{2}*{\scalebox{1.25}{$e^h_{\phi_g}$}}&\multicolumn{1}{|c|}{$\beta=-5\;$}
&4.620E-02  & 1.058E-03 &  5.570E-08 &  3.968E-15\\
&\multicolumn{1}{|c|}{$\beta=\;\;5\;$}
&7.034E-03  & 2.365E-05 &  2.632E-10 &  2.074E-15\\ \hline
\multirow{2}*{\scalebox{1.25}{$e^h_{\varphi_g}$}} &\multicolumn{1}{|c|}{$\beta=-5\;$}
&1.025E-01 &  1.402E-03 &  8.244E-08  & 4.445E-15  \\
&\multicolumn{1}{|c|}{$\beta=\;\;5\;$}
 &1.263E-02 &  3.239E-05 &  3.161E-10  & 1.703E-15\\
\hline\hline
\multicolumn{2}{c|}{GF-DST ($L =8$)}&  $h = 1$  & $h =1/2$ & $h=1/4$ & $h = 1/8$ \\ \hline
\multirow{2}*{\scalebox{1.25}{$e^h_{\phi_g}$}}&\multicolumn{1}{|c|}{$\beta=-5\;$}
 &4.823E-02 &  1.112E-03  & 3.139E-05  & 3.133E-05  \\
&\multicolumn{1}{|c|}{$\beta=\;\;5\;$}
&8.183E-03 &  7.245E-05  & 5.317E-05  & 5.381E-05 \\ \hline
\multirow{2}*{\scalebox{1.25}{$e^h_{\varphi_g}$}}&\multicolumn{1}{|c|}{$\beta=-5\;$}
 & 6.613E-02 &  5.159E-02 &  5.159E-02 &  5.159E-02 \\
&\multicolumn{1}{|c|}{$\beta=\;\;5\;$}
& 6.840E-02 &  6.840E-02 &  6.840E-02 &  6.840E-02 \\
\hline\hline
\multicolumn{2}{c|}{GF-DST ($h =1/8$)}&  $L = 8$  & $L = 16$ &  $L = 32$&  $L = 64$ \\ \hline
\multirow{2}*{\scalebox{1.25}{$e^h_{\phi_g}$}}&\multicolumn{1}{|c|}{$\beta=-5\;$}
& 3.133E-05& 3.848E-06  & 4.789E-07 &  5.980E-08 \\
&\multicolumn{1}{|c|}{$\beta=\;\;5\;$}
& 5.381E-05& 6.212E-06  & 7.606E-07 &  9.445E-08 \\ \hline
\multirow{2}*{\scalebox{1.25}{$e^h_{\varphi_g}$}}&\multicolumn{1}{|c|}{$\beta=-5\;$}
 & 5.159E-02& 2.572E-02  & 1.072E-02 &  5.248E-03 \\
&\multicolumn{1}{|c|}{$\beta=\;\;5\;$}
& 6.840E-02& 3.398E-02  & 1.415E-02 &  6.928E-03\\
\end{tabularx}
{\rule{\temptablewidth}{1pt}}
\end{center}
\end{table}

\begin{table}[htbp]
\tabcolsep 0pt \caption{Different energies of the ground state and $I^h$ for the NLSE with the 2D
Coulomb interaction  for different $\beta$.}
\label{virial_2d_Coulomb}
\begin{center}\vspace{-0.5em}
\def\temptablewidth{1\textwidth}
{\rule{\temptablewidth}{1pt}}
\begin{tabularx}{\temptablewidth}{@{\extracolsep{\fill}}p{1.25cm}lccccc}
$\beta$  &$E_g$ &$\mu_g$ &$E_{\textrm {kin}}^g$  & $E_{\textrm {pot}}^g$ & $E_{\textrm {int}}^g$ &$ I^h $ \\
\hline
$-10$ &0.1367&-1.4536&1.2611   &4.6592E-01  & -1.5903   &  1.89E-10\\
$-5$  &0.8698&0.1933 &9.4226E-01   &6.0401E-01 & -6.7651E-01  &  2.37E-10\\
$-1$  &1.3808&1.2600 &7.8098E-01   &7.2058E-01 & -1.2080E-01  &  2.60E-10\\
$ 1$  &1.6163&1.7311 &7.2201E-01   &7.7942E-01 &  1.1483E-01  & -2.61E-10\\
$5$   &2.0551&2.5801 &6.3379E-01   &8.9629E-01 &  5.2501E-01  & -2.65E-10\\
$10$  &2.5557&3.5132 &5.5977E-01   &1.0385  &  9.5748E-01  & -2.69E-10\\
\end{tabularx}
{\rule{\temptablewidth}{1pt}}
\end{center}
\end{table}

%========================================================================
%			GS-Example 1
%========================================================================

{\bf Example 3.3:} {\sl The NLSE with the Poisson potential  in 2D.} We take $d=2$ and
$U(\bx)= U_{\rm Lap}(\bx)$ in (\ref{sps-eq})-(\ref{sps-colb}).
The ground state is computed numerically on
a bounded domain $\Omega=[-8,8]^2$ with different mesh size $h$.  Table \ref{gs_2d_Poisson} shows
the errors $e_{\phi_g}^h$ and $e_{\varphi_g}^h$ with $V(\bx) = \frac{1}{2}(x^2+4y^2) $
in (\ref{sps-eq}) for different numerical methods, $\beta$ and mesh size $h$.
In addition, Table \ref{virial_2d_Poisson}
lists the energy $E_g:=E(\phi_g^h)$, chemical potential $\mu_g:=\mu(\phi_g^h)$,
kinetic energy $E_{\rm kin}^g:=E_{\rm kin}(\phi_g^h)$, potential energy $E_{\rm pot}^g:=E_{\rm pot}(\phi_g^h)$,
interaction energy $E_{\rm int}^g:=E_{\rm int}(\phi_g^h)$
and $I^h$ with $h=1/8$ and $V(\bx) =\frac{1}{2}(x^2+4y^2)$  in (\ref{sps-eq}) for different $\beta$.

\

From Tables \ref{tab:gs:SPS:3d}-\ref{virial_2d_Poisson} and additional numerical results
not shown here for brevity, we can see that:
(i) The {\sl GF-NUFFT} method is spectrally accurate in space, while the
{\sl GF-DST} method has a saturation accuracy for a fixed domain;
(ii) The saturation error of the {\sl GF-DST} depends inversely on the domain size $L$,
and it can only reach satisfactory accuracy for some large $L$;
(iii) High accuracy, i.e., 9-digit accurate, is achieved by {\sl GF-NUFFT} as quite expected in the energies, which,
in another way, manifest the high-accuracy advantage of our NUFFT solver.

%%%%%%%    NLS with 2D Coulomb interaction  %%%%%%%%
\begin{table}[h!]
\tabcolsep 0pt \caption{Errors of the ground state for the NLSE with the 2D
Poisson potential with mesh size $h$. } \label{gs_2d_Poisson}
\begin{center}\vspace{-1em}
\def\temptablewidth{1\textwidth}
{\rule{\temptablewidth}{1pt}}
\begin{tabularx}{\temptablewidth}{@{\extracolsep{\fill}}p{0.025cm}lllll}
\multicolumn{2}{c|}{GF-NUFFT }  &  $h = 1$  & $h =1/2$ & $h=1/4$ & $h = 1/8$ \\ \hline
\multirow{2}*{\scalebox{1.25}{$e^h_{\phi_g}$}}&\multicolumn{1}{|c|}{$\beta=-5\;$}
 &2.465E-02  & 1.024E-04 &  4.699E-10  & 2.878E-15\\
&\multicolumn{1}{|c|}{$\beta=\;\;5\;$}
&1.191E-02  & 1.593E-05 &  9.793E-12  & 2.726E-15\\ \hline
\multirow{2}*{\scalebox{1.25}{$e^h_{\varphi_g}$}} &\multicolumn{1}{|c|}{$\beta=-5\;$}
&3.737E-02 &  7.634E-05  & 2.896E-10  & 6.347E-14 \\
&\multicolumn{1}{|c|}{$\beta=\;\;5\;$}
&1.033E-02 &  3.282E-06  & 2.682E-12  & 6.247E-14 \\
\end{tabularx}
{\rule{\temptablewidth}{1pt}}
\end{center}
\end{table}

\begin{table}[htbp]
\tabcolsep 0pt \caption{Different energies of the ground state and $I^h$ for the NLSE with the 2D
Poisson potential  for different $\beta$.}
\label{virial_2d_Poisson}
\begin{center}\vspace{-0.5em}
\def\temptablewidth{1\textwidth}
{\rule{\temptablewidth}{1pt}}
\begin{tabularx}{\temptablewidth}{@{\extracolsep{\fill}}p{1.25cm}lccccc}
$\beta$  &$E_g$ &$\mu_g$ &$E_{\textrm {kin}}^g$  & $E_{\textrm {pot}}^g$ & $E_{\textrm {int}}^g$ &$ I^h $ \\
\hline
$-10$ &1.3533&1.1432&9.8061E-01&  5.8272E-01 &-2.1008E-01&  2.44E-10\\
$-5$  &1.4429&1.3691&8.5784E-01&  6.5889E-01 &-7.3819E-02&  2.54E-10\\
$-1$  &1.4913&1.4819&7.7024E-01&  7.3045E-01 &-9.3826E-03&  2.59E-10\\
$ 1$  &1.5073&1.5139&7.3046E-01&  7.7025E-01 & 6.5762E-03& -2.62E-10\\
$5$   &1.5221&1.5260&6.5959E-01&  8.5854E-01 & 3.9516E-03& -2.70E-10\\
$10$  &1.5076&1.4420&5.8770E-01&  9.8559E-01 &-6.5660E-02& -2.81E-10\\
\end{tabularx}
{\rule{\temptablewidth}{1pt}}
\end{center}
\end{table}

\section{For computing the dynamics}\setcounter{equation}{0}
In this section, we present an efficient and accurate
numerical method for computing the dynamics of the NLSE with
the nonlocal interaction potential (\ref{sps-eq})-(\ref{sps-colb}) and the initial data (\ref{sps-ini})
by combining the NUFFT solver for the nonlocal interaction potential evaluation and the time-splitting Fourier pseudospectral
discretization, and compare it with those existing numerical methods.

\subsection{A numerical method via the NUFFT}

From time $t=t_n$ to $t=t_{n+1}$, the NLSE (\ref{sps-eq}) will
be solved in two splitting
steps. One solves first
\be \label{sps-eq65}
i\,\partial_t \psi(\bx,t)= -\frac{1}{2} \Delta\psi(\bx,t),\qquad \bx \in {\mathbb{ R}}^d,\quad t_n\le t\le t_{n+1},
\ee
for the time step of length $\tau$,
followed by solving
\be \label{sps-eq66}
i\,\partial_t \psi(\bx,t)= \left[V(\bx)+ \beta \,\varphi(\bx,t) \right]
\psi(\bx,t),\quad \varphi(\bx,t) = \left(U\ast |\psi|^2\right)(\bx,t), \qquad \bx \in {\mathbb{ R}}^d,\quad t_n\le t\le t_{n+1},
\ee
for the same time step.
For $t\in[t_n,t_{n+1}]$, Eq. (\ref{sps-eq66}) leaves $|\psi|$
invariant in $t$ \cite{BC7,BJNY}, i.e., $|\psi(\bx,t)|=|\psi(\bx,t_n)|$,
and thus $\varphi$ is time invariant, i.e., $\varphi(\bx,t)=\varphi(\bx,t_n):=\varphi^n(\bx)$,
therefore it  becomes
\be \label{sps-eq67}
i\,\partial_t \psi(\bx,t)= \left[V(\bx)+ \beta \,\varphi^n(\bx) \right]
\psi(\bx,t),\quad \varphi^n(\bx) = \left(U\ast |\psi^n|^2\right)(\bx),
\qquad \bx \in {\mathbb{ R}}^d,\quad t_n\le t\le t_{n+1},
\ee
where $\psi^n(\bx):=\psi(\bx,t_n)$,  which immediately implies that
\be\label{sps-eq68}
\psi(\bx,t)=e^{-i\left[V(\bx)+\beta \,\varphi^n(\bx)\right](t-t_n)}\psi(\bx,t_n),
\quad
\qquad \bx \in {\mathbb{ R}}^d,\quad t_n\le t\le t_{n+1}.
\ee
Then an efficient and accurate numerical method can be designed
by:
(i) adopting  a  second-order Strang splitting \cite{Strang}
or a fourth-order time splitting method \cite{Yoshida} to decouple the nonlinearity;
(ii) truncating the problem on a bounded computational domain $\Omega$, and imposing the
periodic BC on $\partial\Omega$ for the subproblem (\ref{sps-eq65});
(iii) discretizing (\ref{sps-eq65}) in space by the Fourier spectral
method and integrating in time {\it exactly};
(iv) evaluating the nonlocal interaction
$\varphi^n(\bx)$ in (\ref{sps-eq68})
by the algorithm via the NUFFT that discussed in previous sections, and
integrating in time {\it exactly} for (\ref{sps-eq68}).
The details are omitted here for brevity and this method is referred as
the {\sl TS-NUFFT} method.
%We remark here that $\widehat{|\psi^n|^2}({\bf 0})
%\equiv\widehat{|\psi_0|^2}({\bf 0})$
%for $n\ge0$.

For comparison, for the nonlocal interaction in 3D/2D,
when the NUFFT in the above method is replaced by the standard FFT,
we refer the method as {\sl TS-FFT}. In addition,
when the nonlocal interaction $\varphi^n(\bx)$ in (\ref{sps-eq68})
is reformulated as its equivalent PDE formulation
(\ref{PoiDiff})-(\ref{SqrtPoiDiff}) on  $\Omega$
with homogeneous Dirichlet
BC  on $\partial\Omega$
and then discretized by the sine pseudospectral method with an evaluation of
(\ref{sps-eq65}) via the sine spectral method and integrated in time
{\it exactly} \cite{DipJCP,SPMCompare}, we refer it as {\sl TS-DST}.

\subsection{Numerical comparisons}

Again, in order to compare the {\sl TS-NUFFT} method with  the {\sl GF-DST} method for
computing the dynamics, we denote $\rho(\bx,t)=|\psi(\bx,t)|^2$ and
$\varphi(\bx,t)=(U\ast|\psi|^2)(\bx,t)$ and introduce
the errors
\beas
e_{\psi}^h(t)&:=&\frac{\max_{\bx\in \Omega^h}|\psi(\bx,t)-\psi^n_{h}(\bx)|}{\max_{\bx\in \Omega^h}|\psi(\bx,t)|},\qquad
 e_{\varphi}^{h}(t):=\frac{\max_{\bx\in \Omega^h} |\varphi(\bx,t)-\varphi^n_{h}(\bx)|}{\max_{\bx\in \Omega^h}|\varphi(\bx,t)|},\\
e_{\rho}^{h}(t)&:=&\frac{\max_{\bx\in \Omega^h} |\rho(\bx,t)-\rho^n_{h}(\bx)|}
{\max_{\bx\in \Omega^h}|\rho(\bx,t)|}, \qquad t=t_n, \qquad n\ge0,
\eeas
where $\psi^n_{h}(\bx)$, $\varphi^n_h(\bx)$ and $\rho^n_{h}(\bx)$ are obtained numerically by a numerical method
as the approximations of  $\psi(\bx,t)$, $\varphi(\bx,t)$ and $\rho(\bx,t)$ at $t=t_n$, respectively
with a given mesh size $h$ and a very small time step $\tau>0$.
The external potential in (\ref{sps-eq}) and the initial data in (\ref{sps-ini}) are chosen as
\be\label{DY-EXP-INI}
V(\bx)=\frac{|\bx|^2}{2},\qquad \psi(\bx,0)=\psi_0(\bx)=e^{-\frac{|\bx|^2}{2}},\qquad\bx\in
{\mathbb R}^d\;\; {\rm with}\;\;d=3\;{\rm or}\;2.
\ee
In the comparisons, the  ``exact" solution $\psi(\bx,t)$  (and thus $\varphi(\bx,t)$ and $\rho(\bx,t)$)
was obtained numerically via the TS-NUFFT method on a  large enough domain $\Omega$ with very
small enough mesh size $h$ and time step $\tau$. In our computations, we use
the fourth-order time-splitting method  for time integration \cite{Yoshida}.

\medskip
%==================================================================
%		Dynamics--3D Coulomb Case
%==================================================================
{\bf Example 4.1:} {\sl The NLSE with the 3D Coulomb interaction.} Here $d=3$ and
$U(\bx)= U_{\rm Cou}(\bx)$ in (\ref{sps-eq})-(\ref{sps-colb}).
The problem is solved numerically on a bounded computational domain $\Omega=[-8,8]^3$
with time step $\tau =10^{-3}$ and different mesh size $h$.
Table \ref{dy_3d_Coulomb}  list the errors of the wave-function, the density and the 3D Coulomb interaction
at $t=1/8$ obtained by the {\sl TS-NUFFT} and {\sl TS-DST} methods
for different mesh size $h$ and
interaction constant $\beta$.

%---------------------------------------------------------------------------------------------------------------------------
%			Beginning of lists of  tables
%---------------------------------------------------------------------------------------------------------------------------

%%%%%%% Dynamics computation of the 3D SPS by TS-NUFFT and TS-DST solver method %%%%%%%%
\begin{table}[h!]
\tabcolsep 0pt \caption{Errors of the wave-function and the nonlocal interaction at $t=1/8$
for the NLSE with the 3D Coulomb interaction. }
\label{dy_3d_Coulomb}
\begin{center}\vspace{-1em}
\def\temptablewidth{1\textwidth}
{\rule{\temptablewidth}{1pt}}
\begin{tabularx}{\temptablewidth}{@{\extracolsep{\fill}}p{0.025cm}lllll}
\multicolumn{2}{c|}{TS-NUFFT} &   $h = 1$  & $h =1/2$ & $h=1/4$  & $h=1/8$ \\ \hline
\multirow{2}*{\scalebox{1.25}{$e^h_{\psi}(1/8)$}}&\multicolumn{1}{|c|}{$\beta=-5\;$}
  &5.461E-03 &  1.011E-05 &  9.297E-12 &  1.492E-13 \\
&\multicolumn{1}{|c|}{$\beta=\;\;5\;$}
 &3.997E-03 &  7.879E-06 &  6.959E-12 &  1.348E-13 \\ \hline
\multirow{2}*{\scalebox{1.25}{$e^h_{\varphi}(1/8)$}} &\multicolumn{1}{|c|}{$\beta=-5\;$}
 &7.890E-03 &  4.466E-06 &  4.745E-12 &  6.992E-14 \\
&\multicolumn{1}{|c|}{$\beta=\;\;5\;$}
 &6.563E-03 &  2.828E-06 &  1.081E-12 &  6.872E-14\\
 \hline \hline

 % -----      Wave-function potential comparison by TS-DST   ----
 \multicolumn{2}{c|}{TS-DST}&  $h = 1$  & $h =1/2$ & $h=1/4$ & $h=1/8$ \\ \hline
\multirow{2}*{\scalebox{1.25}{$e^h_{\psi}(1/8)$}}&\multicolumn{1}{|c|}{$\beta=-5\;$}
&  2.561E-02  & 3.024E-02  & 3.025E-02  & 3.025E-02 \\
&\multicolumn{1}{|c|}{$\beta=\;\;5\;$}
&2.753E-02 &  3.024E-02 &  3.025E-02  & 3.025E-02 \\ \hline

% -----      Density comparison by TS-DST   ----
\multirow{2}*{\scalebox{1.25}{$e^h_{\rho}(1/8)$}}&\multicolumn{1}{|c|}{$\beta=-5\;$}
&5.567E-03  & 1.444E-05 &  2.397E-07  & 2.441E-07\\
&\multicolumn{1}{|c|}{$\beta=\;\;5\;$}
 &5.590E-03  & 1.416E-05 &  2.560E-07  & 2.568E-07 \\ \hline

% -----      Coulomb interaction comparison by TS-DST   ----
\multirow{2}*{\scalebox{1.25}{$e^h_{\varphi}(1/8)$}} &\multicolumn{1}{|c|}{$\beta=-5\;$}
 &1.099E-01  & 1.099E-01  & 1.099E-01  & 1.099E-01\\
&\multicolumn{1}{|c|}{$\beta=\;\;5\;$}
 &1.117E-01 &  1.117E-01 &  1.117E-01  & 1.117E-01 \\

\end{tabularx}
{\rule{\temptablewidth}{1pt}}
\end{center}
\end{table}

\medskip

%==================================================================
%		Dynamics-----2D Coulomb Case
%==================================================================
{\bf Example 4.2:} {\sl The NLSE with the 2D Coulomb interaction.} Here $d=2$ and
$U(\bx)= U_{\rm Cou}(\bx)$ in (\ref{sps-eq})-(\ref{sps-colb}).
The problem is solved numerically on a bounded computational domain $\Omega=[-16,16]^2$
with time step $\tau =10^{-4}$ and different mesh size $h$.
Table \ref{dy_2d_Coulomb} shows the errors of the wave-function and the 2D Coulomb interaction
at $t=0.5$ obtained by the {\sl TS-NUFFT} and {\sl TS-DST} methods
for different mesh size $h$ and interaction constant  $\beta$.

%%%%%%%%%%%    SDM by NUFFT method:   Dynamics---Example 2  %%%%%%%%%%%%%%

\begin{table}[h!]
\tabcolsep 0pt \caption{Errors of the wave-function and the nonlocal interaction at $t=0.5$ for the NLSE with the
2D Coulomb interaction. }
\label{dy_2d_Coulomb}
\begin{center}\vspace{-1em}
\def\temptablewidth{1\textwidth}
{\rule{\temptablewidth}{1pt}}
\begin{tabularx}{\temptablewidth}{@{\extracolsep{\fill}}p{0.025cm}lllll}

\multicolumn{2}{c|}{TS-NUFFT ($L = 16$)} &   $h = 1$  & $h =1/2$ & $h=1/4$  & $h=1/8$ \\ \hline
\multirow{2}*{\scalebox{1.25}{$e^h_{\psi}(0.5)$}}&\multicolumn{1}{|c|}{$\beta=-5\;$}
&1.582E-01 &  7.468E-03 & 4.746E-06  &  2.954E-12 \\
&\multicolumn{1}{|c|}{$\beta=\;\;5\;$}
&5.118E-02 &  7.756E-04 & 2.476E-10  &  1.268E-12\\ \hline
\multirow{2}*{\scalebox{1.25}{$e^h_{\varphi}(0.5)$}} &\multicolumn{1}{|c|}{$\beta=-5\;$}
&2.219E-02 &  4.242E-03  & 4.169E-06  & 3.756E-12\\
&\multicolumn{1}{|c|}{$\beta=\;\;5\;$}
 &3.235E-02 &  2.451E-04  & 3.117E-11  & 7.586E-13\\
 \hline\hline

 % ------- the TS-DST for 2D Coulomb case -----
 \multicolumn{2}{c|}{TS-DST ($L = 16$) } &   $h = 1$  & $h =1/2$ & $h=1/4$  & $h=1/8$ \\ \hline
\multirow{2}*{\scalebox{1.25}{$e^h_{\psi}(0.5)$}}&\multicolumn{1}{|c|}{$\beta=-5\;$}
 &1.175E-01 &  5.576E-02 &  6.311E-02 &  6.312E-02\\
&\multicolumn{1}{|c|}{$\beta=\;\;5\;$}
&6.477E-02 &  6.308E-02 &  6.313E-02 &  6.313E-02\\ \hline
\multirow{2}*{\scalebox{1.25}{$e^h_{\varphi}(0.5)$}} &\multicolumn{1}{|c|}{$\beta=-5\;$}
&4.286E-02 &  2.449E-02  & 2.449E-02 &  2.449E-02\\
&\multicolumn{1}{|c|}{$\beta=\;\;5\;$}
&6.854E-02 &  4.412E-02  & 4.455E-02 &  4.478E-02\\
\hline \hline

 % ------- the TS-DST for 2D Coulomb case with different L  -----
 \multicolumn{2}{c|}{TS-DST ($h = 1/8$)} &  $L= 8$  & $L=16$ & $L=32$ & $L=64$ \\ \hline
\multirow{2}*{\scalebox{1.25}{$e^h_{\psi}(0.5)$}}&\multicolumn{1}{|c|}{$\beta=-5\;$}
& 1.263E-01 &  6.312E-02  & 3.156E-02   &1.578E-02\\
&\multicolumn{1}{|c|}{$\beta=\;\;5\;$}
& 1.264E-01 &  6.313E-02  & 3.156E-02   &1.578E-02\\ \hline
\multirow{2}*{\scalebox{1.25}{$e^h_{\varphi}(0.5)$}} &\multicolumn{1}{|c|}{$\beta=-5\;$}
& 4.907E-02 &  2.449E-02  & 1.021E-02 &  4.999E-03\\
&\multicolumn{1}{|c|}{$\beta=\;\;5\;$}
& 9.038E-02 &  4.500E-02  & 1.875E-02 &  9.181E-03\\

\end{tabularx}
{\rule{\temptablewidth}{1pt}}
\end{center}
\end{table}

\medskip

%==================================================================
%		Dynamics-----2D Poisson Case
%==================================================================
{\bf Example 4.3:} {\sl The NLSE with the 2D Poisson potential.} Here $d=2$ and
$U(\bx)= U_{\rm Lap}(\bx)$ in (\ref{sps-eq})-(\ref{sps-colb}).
Again, the problem is solved numerically on a bounded computational domain $\Omega=[-16,16]^2$
with time step $\tau =10^{-4}$ and different mesh size $h$.
Table \ref{dy_2d_Coulomb} shows the errors of the wave-function and the 2D Coulomb interaction
at $t=0.5$ obtained by the {\sl TS-NUFFT}  method
for different mesh size $h$ and interaction constant  $\beta$.
We remark here that the {\sl TS-DST} method is not applicable for this case \cite{ExtBd,SPMCompare},
therefore here we only present the results for the {\sl TS-NUFFT}  method.

%%%%%%%%%%%    Results on 2D SPS by TS-NUFFT:   Dynamics---Example 1     %%%%%%%%%%%%%%%%%
\begin{table}[h!]
\tabcolsep 0pt \caption{Errors of the wave-function and the Poisson potential at $t=0.5$ for the NLSE with the
2D Poisson potential. }
\label{dy_2d_Poisson}
\begin{center}\vspace{-1em}
\def\temptablewidth{1\textwidth}
{\rule{\temptablewidth}{1pt}}
\begin{tabularx}{\temptablewidth}{@{\extracolsep{\fill}}p{0.025cm}lllll}

\multicolumn{2}{c|}{TS-NUFFT} &   $h = 1$  & $h =1/2$ & $h=1/4$  & $h=1/8$ \\ \hline
\multirow{2}*{\scalebox{1.25}{$e^h_{\psi}(0.5)$}}&\multicolumn{1}{|c|}{$\beta=-5\;$}
&5.833E-02 &  2.599E-04 &  3.211E-09 &  7.524E-13 \\
&\multicolumn{1}{|c|}{$\beta=\;\;5\;$}
&2.658E-02 &  9.083E-05 &  3.395E-12 &  1.124E-12\\ \hline
\multirow{2}*{\scalebox{1.25}{$e^h_{\varphi}(0.5)$}} &\multicolumn{1}{|c|}{$\beta=-5\;$}
 &1.329E-02 &  8.840E-05  & 1.072E-09  & 3.974E-13\\
&\multicolumn{1}{|c|}{$\beta=\;\;5\;$}
&4.645E-03 &  2.805E-06  & 8.322E-13  & 5.821E-13\\

\end{tabularx}
{\rule{\temptablewidth}{1pt}}
\end{center}
\end{table}
%%%%%%%%%%%%   end of 2D SPS by TS-NUFFT:   Dynamics---Example 1  %%%%%%%%%%%

\

From Tables \ref{dy_3d_Coulomb}--\ref{dy_2d_Poisson} and additional numerical results not shown here
for brevity, we can draw the following conclusions:
(i) The {\sl TS-DST}, if applicable, can not resolve the wave-function or the potential very accurately,
while the {\sl TS-NUFFT} achieves the spectral accuracy;
(ii) The saturated accuracy by {\sl TS-DST} decreases as the computation
domain increases;
(iii) As long as for the physical observables, e.g., the density $\rho$,
are concerned, the {\sl TS-DST} method can still capture reasonable accuracy (cf. Tab. \ref{dy_3d_Coulomb}).

%     End of  the Dynamics table list
%---------------------------------------------------------------------------------------------------------------------------

\

\subsection{Applications}

\begin{figure}[t!]
\centerline{
\psfig{figure=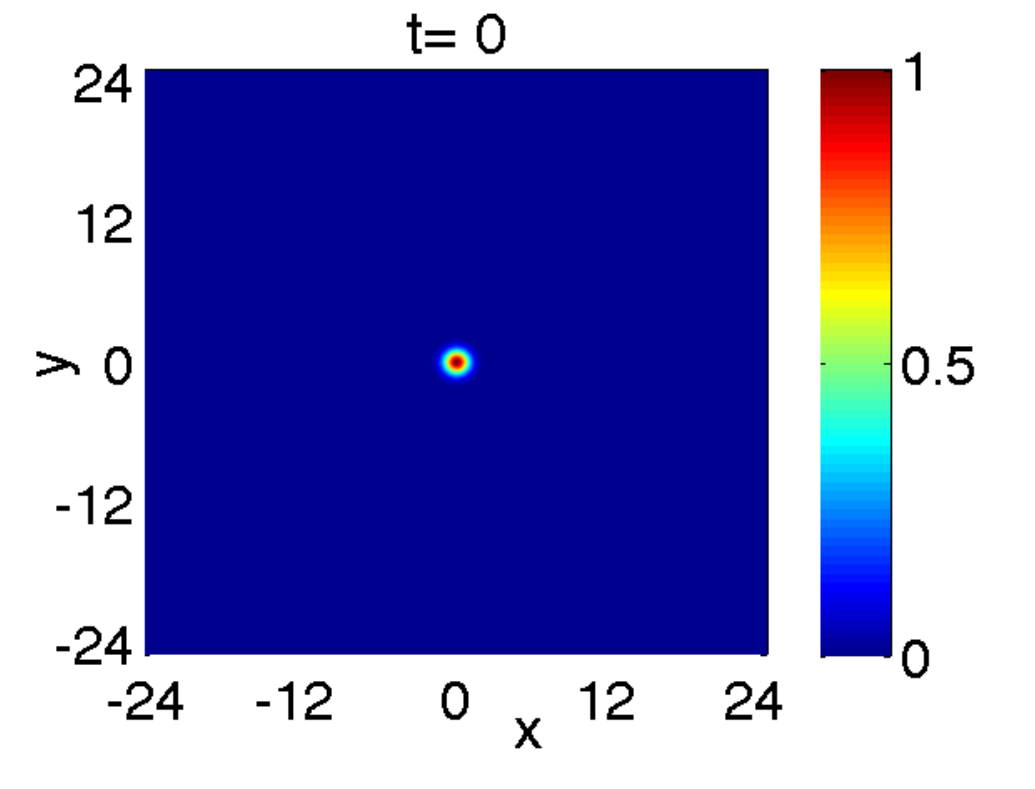,height=4.9cm,width=6.05cm}\qquad
\psfig{figure=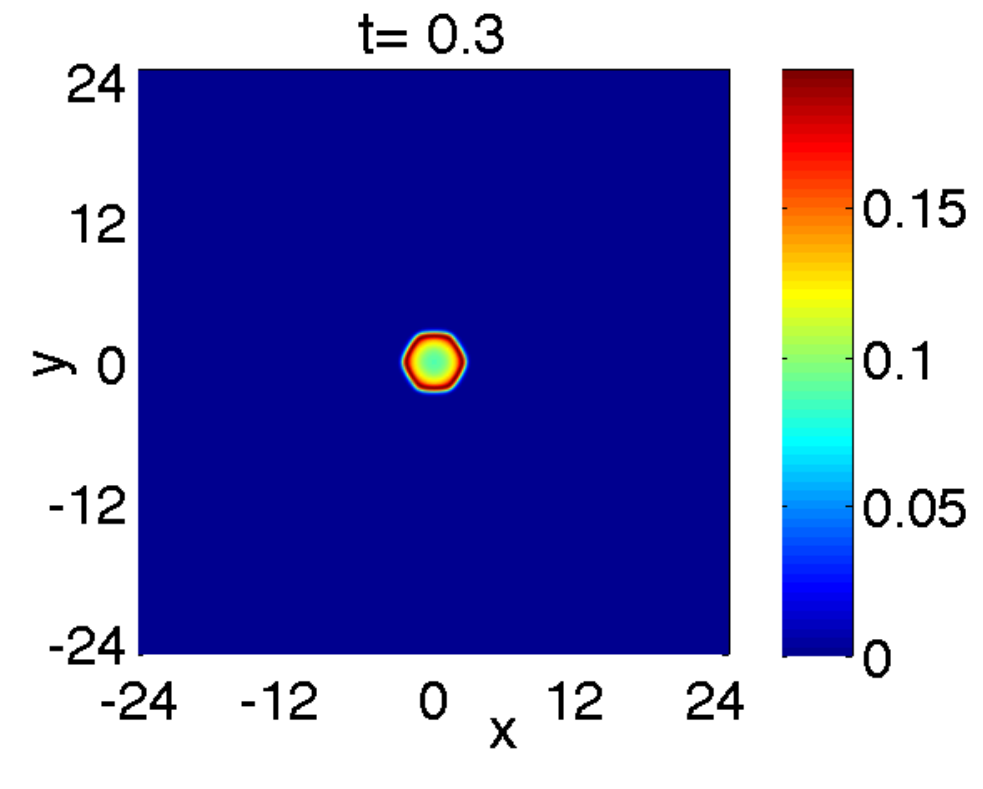,height=4.9cm,width=6.05cm}}
\centerline{
\;\psfig{figure=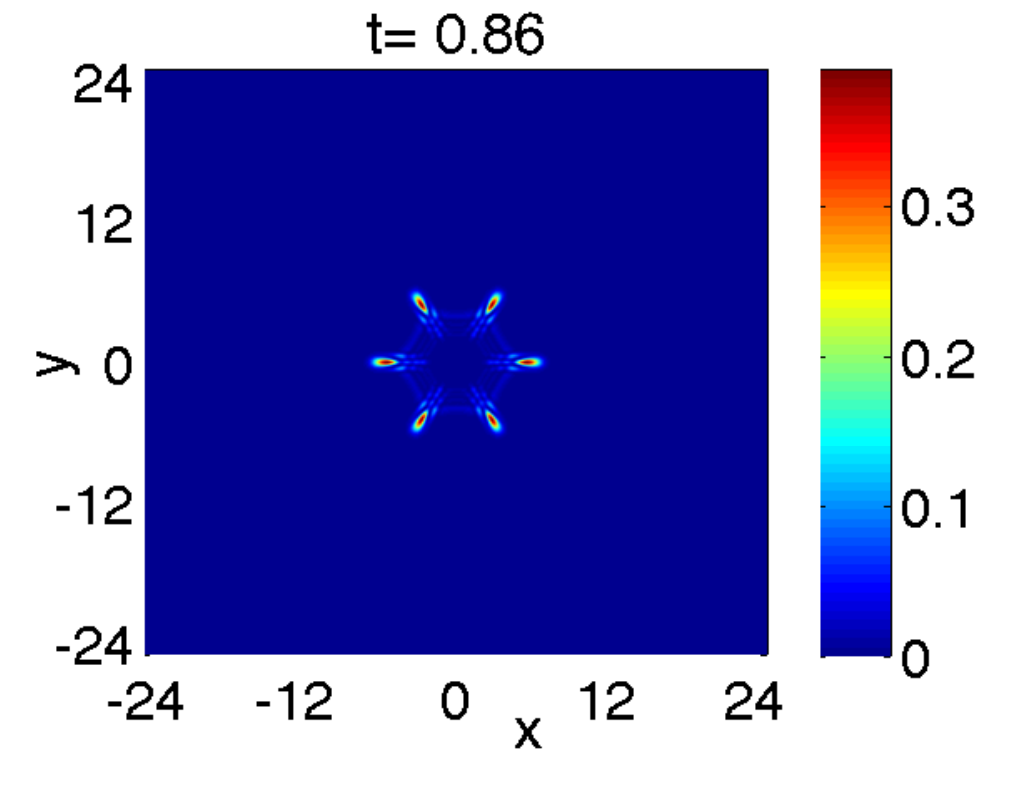,height=4.9cm,width=6.05cm}\qquad
\;\psfig{figure=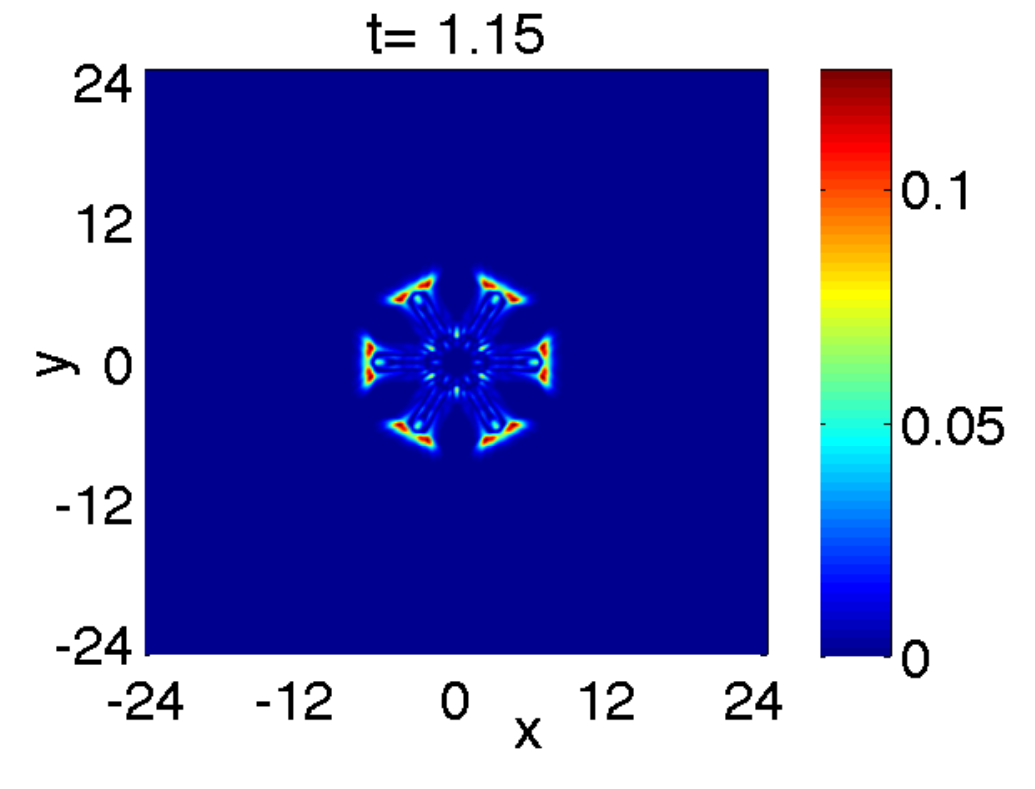,height=4.9cm,width=6.05cm} }
\centerline{
\psfig{figure=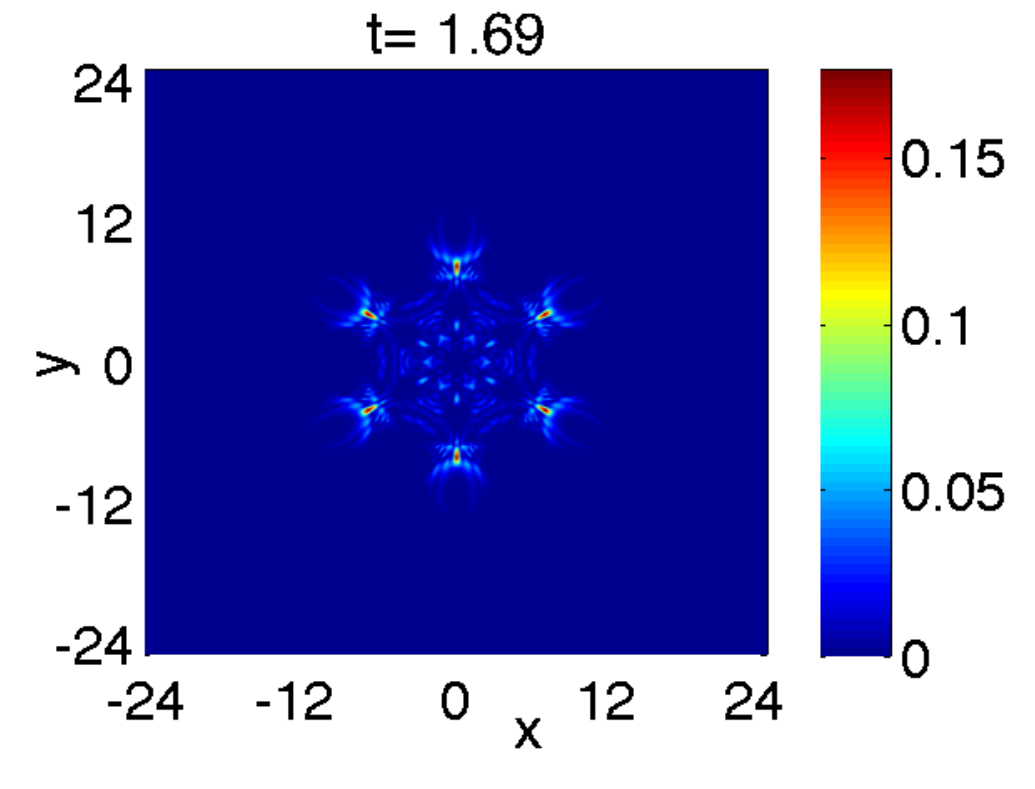,height=4.9cm,width=6.05cm}\qquad
\psfig{figure=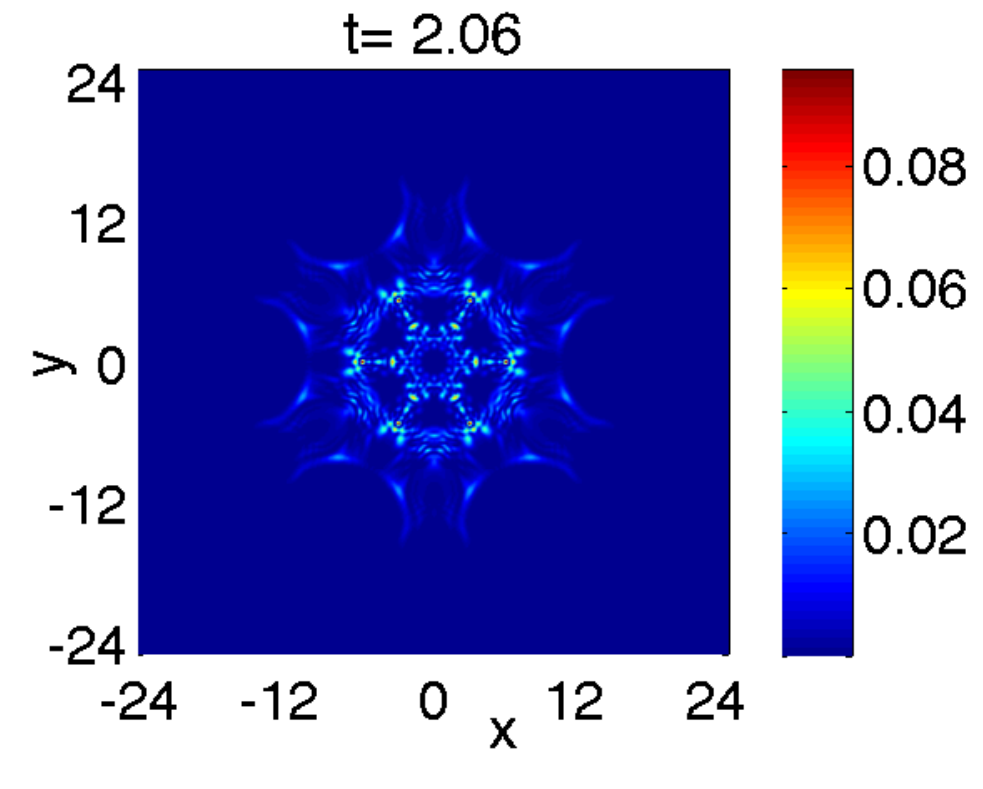,height=4.9cm,width=6.05cm} }
\centerline{
\psfig{figure=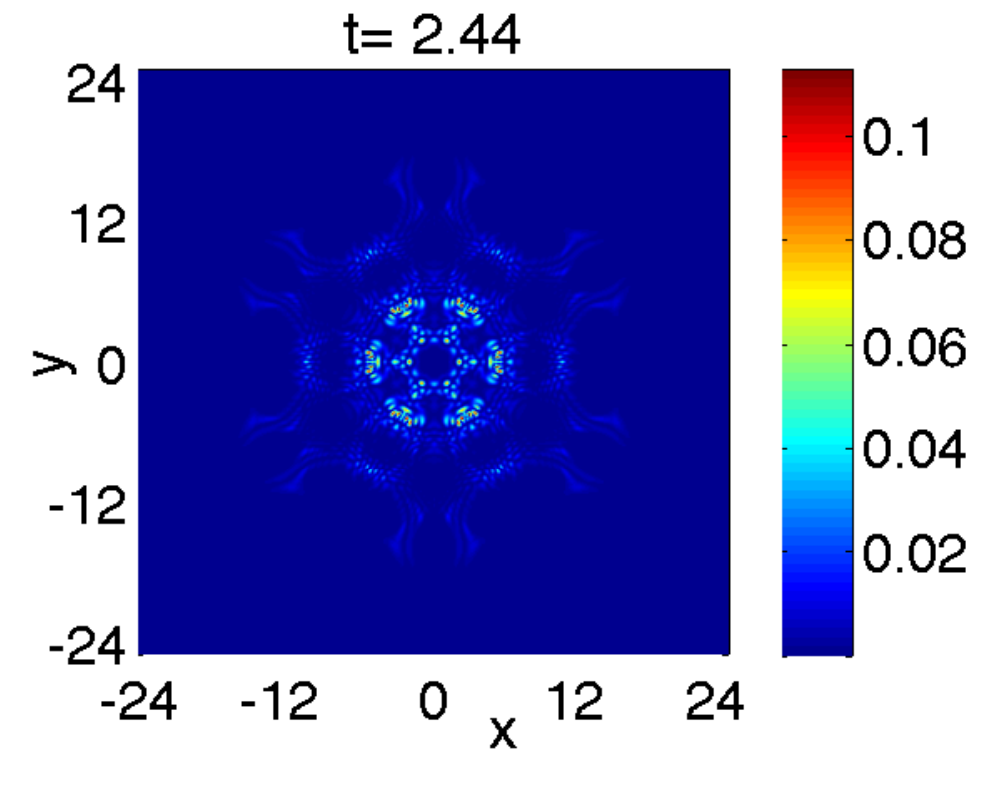,height=4.9cm,width=6.05cm}\qquad
\psfig{figure=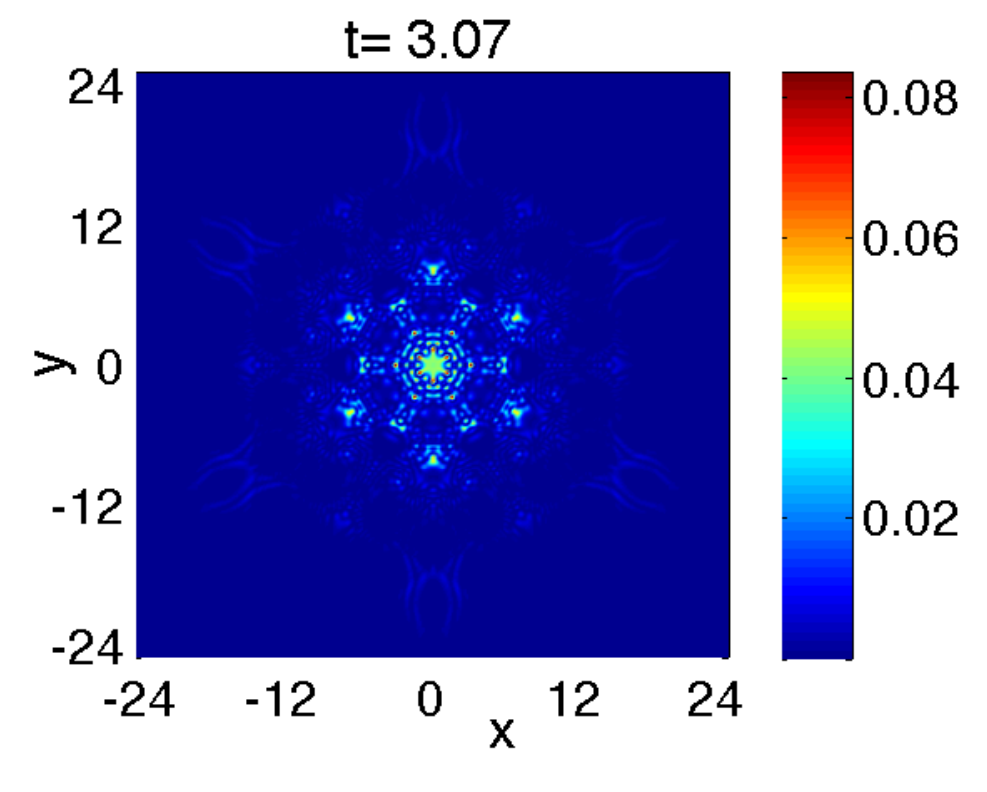,height=4.9cm,width=6.05cm} }
\caption{Contour plots of the density $\rho(x,y,t)$
of the NLSE with the Coulomb interaction and a honeycomb potential
in 2D at different times. }
\label{dy_sdm_app_honey}
\end{figure}

To further demonstrate the efficiency and accuracy of the numerical method via the NUFFT,
we simulate the long-time dynamics of the 2D NLSE with the Coulomb interaction, i.e.,
$d=2$ and  $U(\bx)= U_{\rm Cou}(\bx)$ and $\beta=5$ in (\ref{sps-eq})-(\ref{sps-colb}), and a
honeycomb external potential \cite{BJNY,HoneycombChenWu} defined as
\bea
V(\bx)= 10\left[\cos({\bf b_1}\cdot {\bf x})+\cos({\bf b_2}\cdot {\bf
x})+\cos({\bf(b_1\!+\!b_2)}\cdot {\bf x})\right],\qquad \bx= (x,y)^T\in{\mathbb R}^2,
\eea
with ${\bf b_1} =\frac{\pi}{4}(\sqrt{3},1)^T$ and ${\bf b_2} =\frac{\pi}{4}(-\sqrt{3},1)^T$.
This example can be formally used to describe the dynamics of the electrons in a graphene.
The initial data in  (\ref{sps-ini})  is taken  as $\psi_0(x,y) = e^{-(x^2+y^2)/2}$ for $\bx\in {\mathbb R}^2$
and the problem is solved numerically on $\Omega=[-32,32]^2$ by using the
{\sl TS-NUFFT} with mesh size $h=\frac{1}{16}$ and time step $\tau=10^{-4}$.
Figure \ref{dy_sdm_app_honey} shows the contour plots of the density
$\rho(x,y,t)$ at different times.

\section{Conclusion}
An efficient and accurate numerical method via the NUFFT was proposed for
the fast evaluation of different nonlocal interactions including the Coulomb interactions in 3D/2D
and the interaction kernel taken as either the Green's function of the Laplace operator in 3D/2D/1D
or nonlocal interaction kernels in 2D/1D  obtained from the 3D Schr\"{o}dinger-Poisson system under
strongly external confining potentials via dimension reduction.
The method was compared extensively with those existing numerical methods and was
demonstrated that it can achieve much more accurate numerical results,
especially on a smaller computational domain and/or with anisotropic interaction
density. Eficient and accurate numerical methods were then presented
for computing the ground state and dynamics of the nonlinear Schr\"{o}dinger equation
with nonlocal interactions by combining
the normalized gradient flow with the backward Euler Fourier pseudospectral
discretization and time-splitting Fourier pseudospectral method, respectively,
together with the fast and accurate NUFFT method for evaluating the nonlocal interactions.
Extensive numerical comparisons were carried out between the proposed numerical methods and other existing methods
for studying ground state and dynamics of the NLSE with different nonlocal
interactions. Numerical results showed that the methods via the NUFFT
perform much better than those existing methods in terms of accuracy and efficiency,
especially when the computational domain is chosen smaller and/or
the solution is anisotropic.

\section*{Acknowledgments}
Part of this work was done when the authors were visiting
Beijing Computational Science Research Center in the summer of 2014.
We acknowledge support from  the Ministry of Education of Singapore grant R-146-000-196-112
(W. Bao), the National Science Foundation under grant DMS-1418918 (S. Jiang),
the French ANR-12-MONU-0007-02 BECASIM (Q. Tang)
and  the Austrian Science Foundation (FWF) under
grant No. F41 (project VICOM), grant No. I830 (project LODIQUAS) and the Austrian Ministry of Science
and Research via its grant for the WPI (Q. Tang and Y. Zhang).  The
computation results presented have been achieved by using the Vienna Scientific Cluster.

%%%%%%%%%%%%%%%%%%%%%%%%%%%%%%%%%%%

\end{document}